\documentclass[12pt,oneside,french]{article}
\usepackage{amssymb,amsmath,babel}
\ifx\optionkeymacros\undefined\else \fi

\catcode`\Œ=\active\defŒ{{\aa}}       
\catcode`\º=\active\defº{\int}        
\catcode`\=\active\def{\c c}        
\catcode`\¶=\active\def¶{\partial}    
\catcode`\Ä=\active\defÄ{\oint}       
\catcode`\Æ=\active\defÆ{\triangle}   
\catcode`\Â=\active\defÂ{\neg}        
\catcode`\µ=\active\defµ{\mu}         
\catcode`\¿=\active\def¿{{\o}}        
\catcode`\¹=\active\def¹{\pi}         
\catcode`\Ï=\active\defÏ{{\oe}}       
\catcode`\§=\active\def§{{\ss}}       
\catcode`\ =\active\def {\dagger}     
\catcode`\Ã=\active\defÃ{\sqrt}       
\catcode`\·=\active\def·{\Sigma}      
\catcode`\Å=\active\defÅ{\approx}     
\catcode`\½=\active\def½{\Omega}      
\catcode`\£=\active\def£{{\it\$}}     
\catcode`\°=\active\def°{\infty}      
\catcode`\¤=\active\def¤{{\S}}        
\catcode`\¦=\active\def¦{{\P}}        
\catcode`\¥=\active\def¥{\bullet}     
\catcode`\»=\active\def»{\leavevmode\raise.585ex\hbox{\b a}}      
\catcode`\¼=\active\def¼{\leavevmode\raise.6ex\hbox{\b o}}        
\catcode`\­=\active\def­{\not=}       
\catcode`\²=\active\def²{\leq}        
\catcode`\³=\active\def³{\geq}        
\catcode`\Ö=\active\defÖ{\div}        
\catcode`\É=\active\defÉ{{\dots}}     
\catcode`\¾=\active\def¾{{\ae}}       
\catcode`\Ç=\active\defÇ{\ll}         
\catcode`\Ò=\active\defÒ{``}          
\catcode`\Á=\active\defÁ{!`}          
\catcode`\¢=\active\def¢{\rlap/c}     
\catcode`\Ô=\active\defÔ{`}           
\catcode`\Õ=\active\defÕ{'}           

\catcode`\=\active\def{{\AA}}       
\catcode`\'=\active\def'{\c C}        
\catcode`\¯=\active\def¯{{\O}}        
\catcode`\¸=\active\def¸{\Pi}         
\catcode`\Î=\active\defÎ{{\OE}}       
\catcode`\®=\active\def®{{\AE}}       
\catcode`\×=\active\def×{\diamond}    
\catcode`\¡=\active\def¡{\accent'27}  
\catcode`\Ó=\active\defÓ{''}          
\catcode`\±=\active\def±{\pm}         
\catcode`\È=\active\defÈ{\gg}         
\catcode`\À=\active\defÀ{?`}          
\catcode`\Ð=\active\defÐ{--}          
\catcode`\Ñ=\active\defÑ{---}         

\catcode`\Š=\active\defŠ{\"a}        
\catcode`\'=\active\def'{\"e}        
\catcode`\•=\active\def•{\"{\i}}     
\catcode`\š=\active\defš{\"o}        
\catcode`\Ÿ=\active\defŸ{\"u}        
\catcode`\Ø=\active\defØ{\"y}        
\catcode`\€=\active\def€{\"A}        
\catcode`\…=\active\def…{\"O}        
\catcode`\†=\active\def†{\"U}        
\catcode`\‡=\active\def‡{\'a}        
\catcode`\Ž=\active\defŽ{\'e}        
\catcode`\'=\active\def'{\'{\i}}     
\catcode`\—=\active\def—{\'o}        
\catcode`\œ=\active\defœ{\'u}        
\catcode`\ƒ=\active\defƒ{\'E}        
\catcode`\ˆ=\active\defˆ{\`a}        
\catcode`\=\active\def{\`e}        
\catcode`\"=\active\def"{\`{\i}}     
\catcode`\˜=\active\def˜{\`o}        
\catcode`\=\active\def{\`u}        
\catcode`\Ë=\active\defË{\`A}        
\catcode`\‹=\active\def‹{\~a}        
\catcode`\–=\active\def–{\~n}        
\catcode`\›=\active\def›{\~o}        
\catcode`\Ì=\active\defÌ{\~A}        
\catcode`\"=\active\def"{\~N}        
\catcode`\Í=\active\defÍ{\~O}        
\catcode`\‰=\active\def‰{\^a}        
\catcode`\=\active\def{\^e}        
\catcode`\"=\active\def"{\^{\i}}     
\catcode`\™=\active\def™{\^o}        
\catcode`\ž=\active\defž{\^u}        

\let\optionkeymacros\null

\newtheorem{Theoreme}{Th\'eor\`eme}
\newtheorem{Lemme}{Lemme}
\newtheorem{Proposition}{Proposition}
\newtheorem{Corollaire}{Corollaire}

\newcommand{\NN}{\mathbb N}

\newcommand{\ZZ}{\mathbb Z}
\newcommand{\QQ}{\mathbb Q}
\newcommand{\RR}{\mathbb R}
\newcommand{\CC}{\mathbb C}

\newcommand{\sqm}[4]{\displaystyle{
\left({#1 \atop #3}{#2 \atop #4}\right)}}
\newcommand{\binomial}[2] {{\binom{#1}{#2}}} 
\def\tr{{\bf t}}
\def\no{{\bf n}}

\def\und{\underline}
\def\SL{{\bf SL}}

\def\aa{a}

\begin{document}

\begin{center}

\LARGE{IdŽaux stables dans certains\\
anneaux diffŽrentiels\\ de formes quasi-modulaires de Hilbert.}

\vspace{15pt}

\Large{Federico Pellarin\footnote{Laboratoire LMNO,
Universit\'e de Caen, BP 5186,
14032 Caen Cedex, France.}}

\end{center}

\vspace{15pt}

\noindent {\footnotesize {\bf English abstract.} In \cite{Nesterenko:Modular} Nesterenko 
proved, among other results, the algebraic independence over $\QQ$ of the numbers $\pi,e^{\pi},\Gamma(1/4)$.
A very important feature of his proof is a multiplicity estimate for quasi-modular forms associated to $\SL_2(\ZZ)$ which 
involves profound differential properties of certain non-linear differential systems.

The aim of this article is to begin the study of the analogous properties for Hilbert modular and quasi-modular forms,
especially those which are associated with the number field $\QQ(\sqrt{5})$.
We show that the differential structure of these functions has several analogies with
the differential structure of the quasi-modular forms associated to $\SL_2(\ZZ)$.}

\section{Introduction.} 
Soit ${\cal A}$ un anneau commutatif, et ${\cal D}$ un ensemble de dŽrivations sur ${\cal A}$. Un idŽal ${\cal I}$ de ${\cal
A}$ est {\em stable} pour ${\cal D}$ (ou {\em ${\cal D}$-stable}, ou encore, {\em stable}
si la rŽfŽrence ˆ ${\cal D}$ est claire) si pour tout $D\in{\cal D}$ on a
\[D{\cal I}:=\{Dx,\mbox{ avec }x\in{\cal I}\}\subset{\cal I}.\]

\medskip

\noindent {\bf DŽfinition.} On dit que l'anneau diffŽrentiel $({\cal A},{\cal D})$ satisfait la {\em propriŽtŽ
de Ramanujan}, s'il existe un ŽlŽment non nul $\kappa\in{\cal A}$ tel que tout idŽal premier non nul ${\cal I}$ de ${\cal A}$
qui est ${\cal D}$-stable contient $\kappa$.

\medskip

Soit $q$ un nombre complexe tel que $|q|<1$, soient $E_2(q),E_4(q),E_6(q)$
les sŽries de Fourier des sŽries d'Eisenstein classiques de poids $2,4,6$ (\footnote{NotŽes respectivement $3\pi^{-2}G_1$
p. 156 et $E_2,E_3$ p. 151 de \cite{Serre:Cours}.}) et posons $D=q(d/dq)$.
On sait alors que $D$ induit une dŽrivation sur ${\cal R}=\CC[E_2,E_4,E_6]$. 
On peut montrer que l'anneau diffŽrentiel $({\cal R},D)$ satisfait la
propriŽtŽ de Ramanujan~: tout idŽal premier non nul $D$-stable de ${\cal R}$ contient 
$\Delta=E_4^3-E_6^2$.

En fait, aussi l'anneau ${\cal R}[q]$ est un anneau diffŽrentiel, et on peut dŽmon\-trer qu'il satisfait la
propriŽtŽ de Ramanujan (voir \cite{Pellarin:Bordeaux}).

La propriŽtŽ de Ramanujan est indispensable dans les estimations de multiplicitŽ. Par exemple,
Nesterenko a utilisŽ une propriŽtŽ plus faible que la propriŽtŽ de 
Ramanujan pour $({\cal R}[q],D)$ pour dŽmontrer ce qui suit (cf. \cite{Nesterenko:Modular} et \cite{Nesterenko:Introduction},
chapitres 3, 10).

\medskip

\begin{Theoreme} Il existe une constante explicite $c>0$ telle que
pour tout polyn™me $P\in{\cal R}$ de degrŽ au plus $N$, on ait~:
\[\mbox{ord}_{q=0}P\leq cN^4.\]\label{theoreme:mult_1}\end{Theoreme}

Quant'aux formes modulaires de Hilbert de deux (ou de plusieurs variables complexes),
ce serait Žvidemment trs intŽressant d'Žtablir des
estimations de multiplicitŽ gŽnŽralisant le thŽorme \ref{theoreme:mult_1} pour ces fonctions~; avec ce texte,
nous voudrions inaugurer cette Žtude.

Une premire motivation de ce texte est d'introduire une notion de {\em forme quasi-modulaire de Hilbert} 
gŽnŽralisant celle de \cite{Zagier:Quasimodulaire}~; elle
ne semble pas, ˆ prŽsent, appara"tre dans la literature. 

Soit $K$ un corps de
nombres totalement rŽel de degrŽ $n$ sur $\QQ$, soit $\Gamma_K$ le groupe modulaire de Hilbert associŽ
(cf. \cite{Geer:Hilbert}, p. 5).
Nous allons dŽfinir, au paragraphe \ref{section:quasi_mod}, un anneau ${\cal Y}(K)$
de formes quasi-modulaires de Hilbert pour $\Gamma_K$~; il s'agit de 
fonctions holomorphes de $n$ variables complexes $z_1,\ldots,z_n$. 

Notons~:
\[{\cal D}=\{D_1,\ldots,D_n\}=\{(2\pi{\rm
i})^{-1}\partial/\partial z_1,\ldots, (2\pi{\rm i})^{-1}\partial/\partial z_n\}.\]
Nous verrons que $({\cal Y}(K),{\cal D})$ est un anneau diffŽrentiel.

Nous aurons besoin de quelques propriŽtŽs ŽlŽmentaires de cet anneau, que nous dŽcrirons 
au paragraphe \ref{section:quasi_mod}, et qui seront pour la plupart dŽmontrŽes dans l'appendice.
En particulier, si $n>1$, alors ${\cal Y}(K)$ n'est pas de type fini, contrairement au cas
$K=\QQ$, dŽcrit dans \cite{Zagier:Quasimodulaire}. Voici maintenant nos
rŽsultats concernant les idŽaux premiers ${\cal D}$-stables de ${\cal Y}(K)$.

\begin{Theoreme} Si $n>1$, alors l'anneau ${\cal Y}(K)$ ne possde pas d'idŽaux 
principaux ${\cal D}$-stables non nuls.\label{theoreme:pas_principaux}\end{Theoreme}

Nous avons obtenu des rŽsultats partiels supplŽmentaires sur la structure
des idŽaux premiers diffŽrentiellement stables de ${\cal Y}(K)$ dans le cas 
o $K=\QQ(\sqrt{5})$.

\begin{Theoreme} Supposons que $K=\QQ(\sqrt{5})$~; il existe une forme modulaire de Hilbert non nulle
$\kappa$ avec la propriŽtŽ suivante.
Soit ${\mathfrak P}$ un idŽal premier non nul $D$-stable de 
${\cal Y}(K)$, supposons que ${\mathfrak P}$ contienne une forme modulaire de Hilbert 
non nulle~; alors $\kappa\in{\mathfrak P}$.\label{theoreme:modulaires}\end{Theoreme}

Le thŽorme \ref{theoreme:modulaires} laisse penser que la propriŽtŽ de Ramanujan 
est valide pour ${\cal Y}(\QQ(\sqrt{5}))$.

Voici un autre argument en faveur de cette hypothse. Dans \cite{Pellarin:Bordeaux}, 
nous avons dŽmontrŽ que l'anneau diffŽrentiel ${\cal Y}(\QQ)$ satisfait la propriŽtŽ de Ramanujan.
Tout idŽal premier non nul et diffŽrentiellement stable contient la forme modulaire~:
\[\kappa(\QQ)=\det\sqm{4e_4}{6e_6}{\displaystyle{\frac{de_4}{dz}}}{\displaystyle{\frac{de_6}{dz}}},\] o $e_4,e_6$ 
dŽsignent les sŽries d'Eisenstein de poids $4, 6$ normalisŽes usuelles pour $\SL_2(\ZZ)$.

La forme $\kappa(\QQ)$ est non nulle et proportionnelle ˆ l'unique forme para\-bolique normalisŽe non nulle
de poids $12$ pour $\SL_2(\ZZ)$ (donc proportionnelle ˆ $\Delta(e^{2\pi{\rm i}t})$, $\Im(t)>0$). 

La forme $\kappa(\QQ)$ appara"t donc
comme image d'un {\em crochet de Rankin} de deux gŽnŽrateurs de l'anneau des formes modulaires 
pour $\SL_2(\ZZ)$, et cette expression
joue un r™le important dans la dŽmonstration de cette propriŽtŽ.

Dans le cas o $K=\QQ(\sqrt{5})$, un phŽnomne analogue appara"t.
L'ŽlŽment $\kappa=\kappa(\QQ(\sqrt{5}))$ du thŽorme \ref{theoreme:modulaires}
s'exprime comme un {\em crochet trilinŽaire} (voir \cite{Pellarin:Hilbert}
pour les propriŽtŽs de base de ces crochets) des gŽnŽrateurs $\varphi_2,\chi_5,\chi_6$ de
l'anneau des formes modulaires de Hilbert de {\em poids parallle} (cf. dŽfinition dans ce texte,
ou \cite{Pellarin:Hilbert})~:
\[\kappa(\QQ(\sqrt{5}))=\displaystyle{\det\left(\begin{array}{ccc}
2\varphi_2 & 5\chi_5 & 6\chi_6 \\ & & \\ \displaystyle{\frac{\partial \varphi_2}{\partial z}} & \displaystyle{\frac{\partial
\chi_5}{\partial z}} & 
 \displaystyle{\frac{\partial \chi_6}{\partial z}} \\  & & \\
\displaystyle{\frac{\partial \varphi_2}{\partial z'}} & \displaystyle{\frac{\partial \chi_5}{\partial z'}} & 
 \displaystyle{\frac{\partial \chi_6}{\partial z'}}
\end{array}\right)}.\]
Les crochets multilinŽaires (dont la dŽfinition gŽnŽrale appara"t 
dans \cite{Pellarin:Hilbert}, paragraphe 6.1) mettent en
lumire une analogie entre les formes modulaires
$\kappa(\QQ)$ et $\kappa(\QQ(\sqrt{5}))$~; de plus, dans \cite{Pellarin:Hilbert}, nous avons dŽmontrŽ 
que $\kappa(\QQ(\sqrt{5}))$ est proportionnelle ˆ la forme modulaire $\chi_{15}$
construite par Gundlach dans \cite{Gundlach:Bestimmung1}. Cette analogie a en partie dŽjˆ ŽtŽ remarquŽe par
Resnikoff dans \cite{Resnikoff:Automorphic}, car il avait notŽ que du point de vue diffŽrentiel, 
la sŽrie d'Eisenstein de poids $(2,2)$ normalisŽe $\varphi_2$ pour $\Gamma_{\QQ(\sqrt{5})}$ joue un r™le similaire 
ˆ celui de la sŽrie $E_2$. Dans ce texte, ces analogies seront rendues encore plus impressionantes,
mais leur raison demeure obscure.

Ainsi, les expressions de $\kappa(\QQ)$
et $\kappa(\QQ(\sqrt{5}))$
en termes de crochets multilinŽaires n'expliquent pas compltement pourquoi ces formes 
modulaires jouent un r™le aussi privilŽgiŽ et similaire dans les propriŽtŽs diffŽrentielles
de ${\cal Y}(K)$ pour $K=\QQ$ ou $K=\QQ(\sqrt{5})$. 

En lisant cet article, le lecteur verra qu'en partie, le thŽorme
\ref{theoreme:modulaires} est dž ˆ une co•ncidence que nous ne savons pas
expliquer ˆ prŽsent, et que nous Žsperons d'Žlucider dans des autres travaux.

\subsection{Plan de l'article.}

Voici le plan de cet article~: il y a deux parties de longueur inŽgale, dont les
contenus
sont assez indŽpendants. 
Dans le paragraphe
\ref{section:quasi_mod}, aprs avoir dŽfini l'anneau ${\cal Y}(K)$ et fait une description de ses propriŽtŽs de base, nous
dŽmontrerons le thŽorme
\ref{theoreme:pas_principaux}~; dans la plupart des rŽsultats de ce paragraphe, $K$ est un corps de nombres totalement rŽel de
dŽgrŽ
$n>1$ sur
$\QQ$.

La dŽmonstration du thŽorme \ref{theoreme:modulaires} occupera quant'ˆ elle, tout le reste de 
l'article. Cette dŽmonstration sera ramŽnŽe ˆ la dŽmonstration d'un resultat concernant la structure 
diffŽrentielle d'anneaux de {\em formes modulaires} de Hilbert, ce qui explique la condition
technique du thŽorme \ref{theoreme:modulaires}. Dans ces parties, le choix $K=\QQ(\sqrt{5})$ sera
fixŽ~; les formes quasi-modulaires n'apparaissent pratiquement pas dans cette partie.

\section{Formes quasi-mo\-dulaires.\label{section:quasi_mod}}
Ici, nous dŽfinissons les {\em formes quasi-modulaires de Hilbert} et nous en prŽsentons quelques propriŽtŽs ŽlŽmentaires.

Soit ${\cal H}$ le demi-plan supŽrieur complexe, soit $K$ un corps de nombres totalement
rŽel de degrŽ $n$ sur $\QQ$, soit $\sigma_1,\ldots,\sigma_n$ ses plongements dans $\RR$, soit
${\cal O}_K$ son anneau d'entiers. Si $a\in K$, alors nous Žcrivons aussi $a_i=\sigma_i(a)$.
Le groupe modulaire de Hilbert \[\Gamma_K=\SL_2({\cal O}_K)\] agit sur ${\cal H}^n$ par transformations 
homographiques de la manire usuelle~: soit $z=(z_1,\ldots,z_n)\in{\cal H}^n$ et $\gamma=\sqm{a}{b}{c}{d}\in\SL_2({\cal
O}_K)$. Alors~:
\[\gamma(z)=\left(\frac{a_1z_1+b_1}{c_1z_1+d_1},\ldots,\frac{a_nz_n+b_n}{c_nz_n+d_n}\right).\]
Remarquer que si $n=1$, alors $K=\QQ$ et $\Gamma_K=\SL_2(\ZZ)$.

\medskip

\noindent {\bf Definition.} Soit $F:{\cal H}^n\rightarrow\CC$ une fonction holomorphe.
On dit que $F$ est une {\em forme quasi-modulaire de Hilbert de poids $(k_1,\ldots,k_n)\in\ZZ^n$ et
profondeur $s\in\NN$} si 
\begin{itemize}
\item[(i)] $n=1$ et $F$ est une forme quasi-modulaire pour $\SL_2(\ZZ)$ comme dŽfini par 
Kaneko et Zagier dans
\cite{Zagier:Quasimodulaire}, ou
\item[(ii)] $n>1$, et il existe un polyn™me $P\in\mbox{Hol}({\cal H}^n)[X_1,\ldots,X_n]$ de degrŽ total
$s$ en les indŽterminŽes $X_1,\ldots,X_n$ (dont les coefficients sont des fonctions holomorphes sur ${\cal H}^n$), tel que,
pour tout $z\in{\cal H}^n$ et $\gamma\in\Gamma_K$~:
\begin{equation}
F(\gamma(z)):=\prod_{i=1}^n(c_iz_i+d_i)^{k_i}P\left(\frac{c_1}{c_1z_1+d_1},\ldots,\frac{c_n}{c_nz_n+d_n}\right)
\label{eq:definition_hilbert}.\end{equation}
\end{itemize}
On voit facilement, en utilisant le fait que $F$ admet une expansion en sŽrie de Fourier, que le terme constant de $P$ par
rapport ˆ
$X_1,\ldots,X_n$ est Žgal ˆ
$F$.

La notion de forme quasi-modulaire de Hilbert gŽnŽralise la notion de forme modulaire de Hilbert. 
Par exemple, une forme quasi-modulaire de Hilbert de profondeur $0$ est une forme modulaire.
De plus, si $F$ est une forme quasi-modulaire de poids $(k_1,\ldots,k_n)$ et profondeur $s$, alors 
$F(z,\ldots,z)$ est une forme quasi-modulaire de poids $\sum_ik_i$ et profondeur $ns$ pour $\SL_{2}(\ZZ)$.

Si $F$ est une forme quasi-modulaire de Hilbert de poids $\und{k}:=(k_1,\ldots,k_n)$, alors $D_iF$ est une forme
quasi-modulaire  de Hilbert de poids \[\und{k}+2_i:=(k_1,\ldots,k_{i-1},k_i+2,k_{i+1},\ldots,k_n).\]
Donc $({\cal Y(K)},{\cal D})$ est un anneau diffŽrentiel. 

Des formes quasi-modulaires de Hilbert non constantes et de poids distincts sont $\CC$-linŽairement indŽpendants,
donc les Žspaces vectoriels complexes ${\cal Y}(K)_{(k_1,\ldots,k_n)}$ engendrŽs par les formes quasi-modulaires de Hilbert 
de poids $(k_1,\ldots,k_n)$ sont linŽairement disjoints, filtrŽs par les profondeurs,
et l'anneau ${\cal Y}(K)$ des polyn™mes en des formes quasi-modulaires de Hilbert de tout poids entier et
toute profondeur est multi-graduŽ par les poids~:
\[{\cal Y}(K)=\bigoplus_{(k_1,\ldots,k_n)\in \ZZ^n}{\cal Y}(K)_{(k_1,\ldots,k_n)}.\]

\medskip

\noindent {\bf DŽmonstration du thŽorme \ref{theoreme:pas_principaux}.}
Supposons par l'absurde qu'il existe deux ŽlŽments $A,F\in{\cal Y}(K)$ et une dŽrivation
$D_i\in{\cal D}$ avec~:
\[D_iF=AF.\]
Nous pouvons supposer que $F$ soit une fonction non constante, ce qui implique $A\not=0$. Il existe des formes 
quasi-modulaires de Hilbert $F_{\und{k}}$ non nulles, de poids distincts $\und{k}$ et $A_{\und{h}}$ de poids
distincts
$\und{h}$, de telle sorte que~:
\[F=\sum_{\und{k}}F_{\und{k}},\quad A=\sum_{\und{h}}A_{\und{h}},\] les sommes Žtant finies. Dans l'ensemble $\{F_{\und{k}}\}$,
choisissons un ŽlŽment non constant $F_{\und{k}_0}$, tel que le poids 
$\und{k}_0=(k_{0,1},\ldots,k_{0,n})$ satisfasse la propriŽtŽ suivante~:
le nombre $v(\und{k}_0):=\sum_{j\not=i}k_{0,j}$ est le plus petit possible, et
$k_{0,i}$ est le plus grand possible pour ce choix de $v(\und{k}_0)$.

On en dŽduit que l'ensemble $\{A_{\und{h}}\}$ contient un ŽlŽment non nul $A_{2_i}$, de poids
\[2_i:={(\underbrace{0,\ldots,0}_{i-1\mbox{ fois
}},2,0,\ldots,0)},\] 
tel que:
\[D_iF_{\und{k}_0}=A_{2_i}F_{\und{k}_0}.\] Ceci constitue une contradiction avec le lemme \ref{lemma:non_existence}
de l'appendice.

Dans l'appendice nous dŽmontrons Žgalement le rŽsultat suivant, qui est aussi une consŽquence 
simple du lemme \ref{lemma:non_existence}.
\begin{Theoreme}
Si $[K:\QQ]>1$, alors l'anneau ${\cal Y}(K)$ des formes quasi-modu\-laires de Hilbert pour 
$\Gamma_K$ n'est pas de type fini.
\label{theorem:theorem1}\end{Theoreme}
Ce thŽorme pour $n=2$ est aussi un corollaire 
d'un rŽsultat de N. Oled Aaiez concernant la structure diffŽrentielle des formes quasi-modulaires
non holomorphes, associŽes ˆ des cycles de Hirzebruch-Zagier compacts sur ${\cal H}^2$ (travail ˆ para"tre).

\section{PropriŽtŽ de Ramanujan de formes modulaires.}

Dans les paragraphes qui suivent, nous munissons un anneau de formes modulaires de Hilbert de deux variables complexes
(associŽ au corps $K=\QQ(\sqrt{5})$) de plusieurs dŽrivations liŽes aux crochets de Rankin.
L'anneau diffŽrentiel ainsi obtenu satisfait lui aussi la propriŽtŽ de Ramanujan~;
ceci impliquera le thŽorme \ref{theoreme:modulaires}.

Faisons tout de suite un survol un peu plus prŽcis de ce que nous allons
montrer. Soit ${\cal O}_K$ l'anneau des entiers de $K$.
Nous allons dŽfinir, sur l'anneau ${\cal
T}$ engendrŽ par les formes modulaires de Hilbert de poids parallles pour le groupe $\Gamma=\Gamma_K$, un ensemble de
dŽrivations \[{\mathfrak D}=\{d_1,d_2,e_1,e_2,f_1,f_2\}.\]
Nous verrons que $({\cal T},{\mathfrak D})$ satisfait
la propriŽtŽ de Ramanujan, mais que cet anneau n'a pas d'idŽaux premiers principaux non nuls et stables, contrairement ˆ
l'anneau $({\cal R},D)$ (il y a \og trop de dŽrivations\fg ).

Nous verrons qu'en dŽfinissant \[{\mathfrak D}^*=\{(d_1+d_2)/2,(e_1+e_2)/2,(f_1+f_2)/2\},\]
l'anneau diffŽrentiel $({\cal T},{\mathfrak D}^*)$ satisfait lui aussi la propriŽtŽ de Ramanujan, et a de plus
un, et un seul idŽal premier principal non nul ${\mathfrak D}^*$-stable (il y a le \og bon nombre de dŽrivations\fg ). 

En revanche, l'anneau diffŽrentiel $({\cal T},\{e_1+e_2\})$, possde deux idŽaux premiers non nuls stables $(\chi_5)$ et
$(\chi_{15})$, mais nous ne savons pas dŽmontrer que cet anneau diffŽrentiel satisfait la propriŽtŽ de Ramanujan (Il n'y a
pas \og assez de dŽrivations\fg ).

Voici les rŽsultats principaux de cette partie, que nous prŽferons enoncer avant mme de dŽcrire les objets utilisŽs.

\begin{Theoreme} Si ${\cal Q}$ est un idŽal premier non nul de ${\cal T}$ qui est ${\mathfrak D}$-stable, alors $\chi_{15}\in{\cal
Q}$.
\label{theoreme:ideaux_stables12}\end{Theoreme}

Nous allons montrer qu'il existe un sous-anneau ${\cal T}^*$ de ${\cal T}$, isomorphe ˆ l'anneau de 
polyn™mes $\CC[X_1,X_2,X_3]$ (pour
des indŽterminŽes $X_1,X_2,X_3$), tel que $({\cal T}^*,{\mathfrak D}^*)$ est un anneau diffŽrentiel. 
On verra que $\chi_{15}\not\in{\cal T}^*$, mais que \[\chi=\chi_{15}^2\in{\cal T}^*,\] et 
que ce polyn™me est irrŽductible.

\begin{Theoreme} Si ${\cal P}$ est un idŽal premier non nul de ${\cal T}^*$ qui est ${\mathfrak D}^*$-stable, alors $ \chi\in{\cal
P}$.
\label{theoreme:ideaux_stables}\end{Theoreme}

Nous donnerons une description plus dŽtaillŽe de ces anneaux diffŽrentiels~; pour l'instant, il est utile de savoir que
le thŽorme \ref{theoreme:ideaux_stables12} s'obtient directement comme corollaire du thŽorme
\ref{theoreme:ideaux_stables}, et que le thŽorme \ref{theoreme:modulaires} n'est qu'un corollaire du thŽorme
\ref{theoreme:ideaux_stables12}.

Des calculs numŽriques plus approfondis nous ont permis de prŽciser les thŽormes 
\ref{theoreme:ideaux_stables12} et \ref{theoreme:ideaux_stables}~; ceci sera abordŽ dans la proposition
\ref{theoreme:ideaux_stables_ht_2}
du paragraphe 
\ref{section:paragrapheB}. Les dŽmonstrations que nous donnerons mettent en relief des analogies
formelles entre l'anneau diffŽrentiel $({\cal R},D)$ et l'anneau diffŽrentiel $({\cal T}^*,{\mathfrak D}^*)$.

\subsection{Structure et dŽrivations modulaires.}

Pour les notions ŽlŽmentaires concernant les formes modulaires de Hilbert associŽes au corps $K=\QQ(\sqrt{5})$, nous 
utiliserons essentiellement \cite{Pellarin:Hilbert}~: nous en reprenons en grande partie les notations
et la terminologie. Nous designons souvent par $\und{r}$ un couple de nombres complexes $(r_1,r_2)$. Si $n$ est un entier,
nous posons $n_1=(n,0)$ et $n_2=(0,n)$~: nous Žcrirons par exemple $2_1=(2,0)$.

Soit $M_{\und{r}}(\Gamma)$ l'espace vectoriel des formes modulaires de Hilbert 
de poids $\und{r}$.
Nous avons~:
\[{\cal T}=\bigoplus_{n\in\NN}M_{(n,n)}(\Gamma).\]
Posons aussi~:
\[{\cal L}=\bigoplus_{\und{r}\in\NN^2}M_{\und{r}}(\Gamma).\]
${\cal L}$ est l'anneau des polyn™mes en des formes modulaires de Hilbert de tout poids (non nŽcessairement parallle)~:
contrairement ˆ ${\cal T}$, cet anneau n'est pas de type fini (cf. \cite{Pellarin:Hilbert}, lemme 16).

Dans la suite, nous considŽrons le poids d'une forme modulaire de poids parallle comme un entier, plut™t que
comme un couple d'entiers.

Nous dŽsignons par ${\cal T}_{{\tiny \mbox{sym}}}$ le sous-anneau de ${\cal T}$ dont les ŽlŽments
sont les formes modulaires $F$ telles que $F\circ{\cal C}=F$, o ${\cal C}$ est 
la symŽtrie ${\cal C}(z,z')=(z',z)$~: c'est le sous-anneau des formes modulaires {\em symŽtriques}.
Une forme modulaire $F$ est dite {\em antisymŽtrique} si elle satisfait $F\circ{\cal C}=-F$.

\medskip

Sur les espaces vectoriels de formes modulaires de Hilbert $M_{\und{r}}$
sont dŽfinis certains opŽrateurs diffŽrentiels, qui sont des 
avatars des {\em crochets} de Rankin et Rankin-Cohen.
Soit $F$ une forme modulaire de Hilbert de poids $\und{f}=(f_1,f_2)$. 
Nous utiliserons les opŽrateurs non linŽaires suivants~:
\begin{eqnarray*} 
(2\pi i)^2\Pi_i F & := & f_iF\frac{\partial^2 F}
{\partial z_i^2}-(f_i+1)\left(\frac{\partial F}{\partial
z_i}\right)^2,\\ (2\pi i)^2\Lambda F & := & F\frac{\partial^2 F}{\partial z_1\partial z_2}-
\frac{\partial F}{\partial
z_1}\frac{\partial F}{\partial z_2}.
\end{eqnarray*}
Un calcul direct sur les facteurs d'automorphie, et la dŽtermination du terme constant dans
les expansions en sŽrie de Fourier, montrent que~:
\begin{eqnarray*}
\Pi_i:M_{\und{f}}(\Gamma) & \rightarrow & S_{2\und{f}+4_i}(\Gamma),\\
\Lambda:M_{\und{f}}(\Gamma) & \rightarrow & S_{2\und{f}+(2,2)}(\Gamma),
\end{eqnarray*}
o $S_{\und{r}}(\Gamma)$ est l'espace vectoriel des formes paraboliques de poids $\und{r}$.

Soit $G$ une forme modulaire de Hilbert de poids $\und{g}=(g_1,g_2)$~: nous posons~:
\[(2\pi i)[F,G]_{1_i}:=  f_iF\frac{\partial G}{\partial
z_i}-g_iG\frac{\partial F}{\partial z_i}=\det\sqm{f_iF}{g_iG}{\displaystyle{\frac{\partial F}{\partial z_i}}}{
\displaystyle{\frac{\partial
G}{\partial z_i}}}.\] Nous avons un opŽrateur bilinŽaire~:
\begin{eqnarray*}
[\cdot,\cdot]_{1_i}:M_{\und{f}}(\Gamma)\times M_{\und{g}}(\Gamma) & \rightarrow &
S_{\und{f}+\und{g}+2_i}(\Gamma).\end{eqnarray*}
Soient $F,G,H$ trois formes modulaires de Hilbert de poids parallles $f,g,h$ respectivement. L'opŽrateur~:
\[(2\pi{\rm i})^2{[}F,G,H{]} =\frac{(g+h)}{2}\det\left(\begin{array}{ccc}
fF & gG & hH \\ & & \\ \displaystyle{\frac{\partial F}{\partial z}} & \displaystyle{\frac{\partial G}{\partial z}} & 
\displaystyle{\frac{\partial H}{\partial z}} \\  & & \\
\displaystyle{\frac{\partial F}{\partial z'}} & \displaystyle{\frac{\partial G}{\partial z'}} & 
\displaystyle{\frac{\partial H}{\partial z'}}\end{array}\right)\]
est trilinŽaire et nous avons~:
\[[\cdot,\cdot,\cdot] :M_{f}(\Gamma)\times M_{g}(\Gamma)\times M_{h}(\Gamma)\rightarrow
S_{f+g+h+2}(\Gamma).\] 
On a aussi l'identitŽ~:
\begin{equation}{[}F,G,H{]} =\frac{1}{2}([F,[G,H]_{1_2}]_{1_1}-[F,[G,H]_{1_1}]_{1_2}),\label{eq:triplo}
\end{equation}
voir \cite{Pellarin:Hilbert}, proposition 3.
Un fait trs important ˆ remarquer est que ces crochets {\em ne sont pas bien dŽfinis} sur des polyn™mes 
quelconque en des formes modulaires, mais seulement sur les formes modulaires, puisque ils dŽpendent du poids.
On a le rŽsultat suivant (cf. thŽorme 6.1 de
\cite{Pellarin:Hilbert}, voir aussi \cite{Hirzebruch:Klein})~:

\begin{Proposition} Soit $\varphi_2$ la sŽrie d'Eisenstein de poids $2$ pour $\Gamma$, normalisŽe 
pour que le terme constant de son dŽveloppement de Fourier ˆ l'infini soit Žgal ˆ $1$. Il existe trois formes 
paraboliques non nulles $\chi_5,\chi_6,\chi_{15}$ de poids parallles respectivement $5,6,15$,
avec $\chi_6,\chi_{15}$ (et $\varphi_2$) symŽtriques et $\chi_5$ antisymŽtrique,
satisfaisant les propriŽtŽs suivantes.
\begin{eqnarray*}
\chi_6 & = & \frac{1}{24}\Lambda\varphi_2,\\
\chi_5^2 & = & \frac{-1}{2880}\varphi_2^{-1}((\Pi_1\varphi_2)(\Pi_2\varphi_2)-9(\Lambda\varphi_2)),\\
\chi_{15} & = & \frac{\sqrt{5}}{22}[\chi_6,\varphi_2,\chi_5] .
\end{eqnarray*}
L'anneau des formes modulaires de Hilbert pour le groupe $\Gamma$,
est Žgal ˆ l'anneau (graduŽ par les poids) des polyn™mes isobares en $\varphi_2,\chi_5,\chi_6,\chi_{15}$,
divisŽ par l'idŽal engendrŽ par la rŽlation~:
\begin{equation}\chi_{15}^2=\chi,\label{eq:klein}\end{equation}
avec
\begin{eqnarray*}
\chi & = & \lambda(50000\chi_5^6-1000\varphi_2^2\chi_6
\chi_5^4+\varphi_2^5\chi_5^4-2\varphi_2^4\chi_6^2\chi_5^2+\\
& &+1800\varphi_2\chi_6^3\chi_5^2+\varphi_2^3\chi_6^4-864\chi_6^5)
\end{eqnarray*} et $\lambda=484/49$.\label{proposition:structure2}\end{Proposition}

Sur ${\cal T}$ on a des involutions $\iota,\varsigma$, dŽfinies par~:
\[\varsigma(\chi_{5})=-\chi_{5},\quad\iota(\chi_{15})=-\chi_{15}.\] 
Nous Žcrirons souvent $\overline{x}:=\iota(x)$~: l'involution
$\iota$ est ${\cal T}^*$-linŽaire (o ${\cal T}^*=\CC[\varphi_2,\chi_5,\chi_6]$) et associe ˆ
une forme modulaire de Hilbert $F(z,z')$ de poids parallle, la forme modulaire $F(\epsilon^2 z',\epsilon'{}^2z)$,
o \[\epsilon=\frac{1+\sqrt{5}}{2}.\]
L'involution
$\varsigma$ est ${\cal T}_{{\tiny \mbox{sym}}}$-linŽaire,
et associe ˆ une forme modulaire de Hilbert $F(z,z')$ de poids parallle, la forme modulaire $F(z',z)$, composŽe de $F$ avec
la symŽtrie ${\cal C}$ (noter que ${\cal T}_{{\tiny \mbox{sym}}}=
\CC[\varphi_2,\chi_5^2,\chi_6,\chi_{15}]$, d'aprs la proposition \ref{proposition:structure2}).

\medskip

\noindent {\bf Remarque.} 
nous avons normalisŽ $\varphi_2,\chi_5,\chi_6,\chi_{15}$ de telle sorte que tous les
coefficients de Fourier de leurs expansions ˆ la pointe ˆ l'infini de la surface modulaire 
de Hilbert associŽe ˆ $K$ soient des entiers rationnels 
premiers entre eux. Ce choix n'est pas celui de \cite{Hirzebruch:Klein} ou \cite{Geer:Hilbert}.
C'est pourquoi la relation $\chi_{15}^2-\chi$ avec $\chi$ comme dans
(\ref{eq:klein}), ne ressemble pas ˆ la formule (16) p. 109 de \cite{Hirzebruch:Klein}
ou ˆ la formule du corollaire 2.2 p. 191 de \cite{Geer:Hilbert}.

\subsection{Les six dŽrivations.}
Nous prŽsentons ici les $6$ dŽrivations $d_1,d_2,e_1,e_2,f_1,f_2$.
Soient $G,H$ deux formes modulaires dans ${\cal T}$ de poids $g,h$. Les applications~:
\begin{eqnarray*}
M_{f}(\Gamma) & \rightarrow & S_{f+g+h+2}\\
F & \mapsto & [F,[G,H]_{1_1}]_{1_2}\\
F & \mapsto & [F,[G,H]_{1_2}]_{1_1}
\end{eqnarray*}
sont dŽfinies sur les ŽlŽments de ${\cal T}$ qui sont homognes
pour la graduation de la proposition \ref{proposition:structure2} (des ŽlŽments isobares, c'est-ˆ-dire des formes modulaires de 
poids parallle). 

ConsidŽrons le cas o $G,H\in\{\varphi_2,\chi_5,\chi_6\}$. Pour tout $r$, nous avons six applications 
linŽaires sur $M_r(\Gamma)$ ainsi dŽfinies~:
\begin{eqnarray*}
d_1(X) = [X,[\varphi_{2},\chi_{5}]_{1_2}]_{1_1} & & d_2(X) = [X,[\varphi_{2},\chi_{5}]_{1_1}]_{1_2}\\
e_1(X) = [X,[\varphi_{2},\chi_{6}]_{1_2}]_{1_1} & & e_2(X) = [X,[\varphi_{2},\chi_{6}]_{1_1}]_{1_2}\\
f_1(X) = [X,[\chi_{5},\chi_{6}]_{1_2}]_{1_1} & & f_2(X) = [X,[\chi_{5},\chi_{6}]_{1_1}]_{1_2}.
\end{eqnarray*} Ces applications sont {\em isobares} de poids parallles $9,10,13$, ou en d'autres termes~:
\begin{eqnarray*}
d_1,d_2:M_r(\Gamma) & \rightarrow & S_{r+9}(\Gamma),\\
e_1,e_2:M_r(\Gamma) & \rightarrow & S_{r+10}(\Gamma),\\
f_1,f_2:M_r(\Gamma) & \rightarrow & S_{r+13}(\Gamma).
\end{eqnarray*}
Posons \[d^*=(d_1+d_2)/2,\quad e^*=(e_1+e_2)/2,\quad f^*=(f_1+f_2)/2,\] et
\[d_*=(d_1-d_2)/2,\quad e_*=(e_1-e_2)/2,\quad f_*=(f_1-f_2)/2,\] de telle sorte que
$d_1=d^*+d_*,d_2=d^*-d_*$, etc. Posons aussi \[{\mathfrak D}^*=\{d^*,e^*,f^*\},\quad{\mathfrak D}_*=\{d_*,e_*,f_*\}.\]

Si $F,G\in{\cal T}$ sont deux formes modulaires de poids $f,g$, alors
$d_1(FG)=Fd_1(G)+Gd_1(F)$. Ainsi l'application linŽaire $d_1$
se prolonge par linŽaritŽ en une dŽrivation de l'anneau ${\cal T}$~; si
$F=\sum_{l=0}^kF_l$ est un ŽlŽment de ${\cal T}$ avec $F_l$ forme modulaire de Hilbert de poids parallle $l$, alors 
nous posons~:
\[d_1(F)=\sum_{l=1}^kd_1(F_l).\]
Par le mme argument, tous les ŽlŽments de ${\mathfrak D},{\mathfrak D}^*,{\mathfrak D}_*$ se prolongent en
des dŽrivations sur l'anneau ${\cal T}$, compatibles avec la graduation par les poids des formes modulaires.

\medskip

Soit ${\cal T}_0=\CC[\varphi_2,\chi_5^2,\chi_6]$.
La proposition \ref{proposition:structure2} dit que ${\cal T}$ a une structure de ${\cal T}_0$-module 
libre de rang $4$. Plus prŽcisement
on a \[{\cal T}={\cal T}_0\oplus\chi_5{\cal T}_0\oplus\chi_{15}{\cal T}_0\oplus\chi_5\chi_{15}{\cal T}_0.\] 
En particulier,
${\cal T}_{{\tiny sym}}={\cal T}_0\oplus\chi_{15}{\cal T}_0$ et ${\cal T}^*={\cal T}_0\oplus\chi_5{\cal T}_0$. 
Dans la proposition
qui suit, nous nous servons de cette ${\cal T}_0$-structure pour dŽterminer les images des sous-modules ${\cal T}_{{\tiny sym}}$
et ${\cal T}^*$ et des autres, par les opŽrateurs de dŽrivation de ${\mathfrak D}^*$ ou de ${\mathfrak D}_*$.
\begin{Proposition} 
On a les inclusions suivantes.
\begin{eqnarray}
d^*({\cal T}_{{\tiny sym}}),e^*(\chi_5{\cal T}_{{\tiny sym}}),f^*({\cal T}_{{\tiny sym}})& \subset & \chi_5
{\cal T}_{{\tiny sym}}\nonumber\\
d^*(\chi_5{\cal T}_{{\tiny sym}}),e^*({\cal T}_{{\tiny sym}}),f^*(\chi_5{\cal T}_{{\tiny sym}})& \subset & 
{\cal T}_{{\tiny sym}}\nonumber\\
d_*({\cal T}_{{\tiny sym}}),e_*(\chi_5{\cal T}_{{\tiny sym}}),f_*({\cal T}_{{\tiny sym}}) & \subset & 
{\cal T}_{{\tiny sym}}\nonumber\\
d_*(\chi_5{\cal T}_{{\tiny sym}}),e_*({\cal T}_{{\tiny sym}}),f_*(\chi_5{\cal T}_{{\tiny sym}})& \subset & 
\chi_5{\cal T}_{{\tiny sym}}\nonumber\\
d^*({\cal T}^*),e^*({\cal T}^*),f^*({\cal T}^*)& \subset & {\cal T}^*\nonumber\\
d^*(\chi_{15}{\cal T}^*),e^*(\chi_{15}{\cal T}^*),f^*(\chi_{15}{\cal T}^*)& \subset & 
\chi_{15}{\cal T}^*\label{eq:ideal_principal_stable}\\
d_*({\cal T}^*),e_*({\cal T}^*),f_*({\cal T}^*)& \subset &\chi_{15}{\cal T}^*\nonumber\\
d_*(\chi_{15}{\cal T}^*),e_*(\chi_{15}{\cal T}^*),f_*(\chi_{15}{\cal T}^*)& \subset & {\cal T}^*\nonumber.
\end{eqnarray}
De plus, on a les ŽgalitŽs~:
\begin{equation}d_2(x)=\overline{d_1(\overline{x})},e_2(x)=\overline{e_1(\overline{x})},f_2(x)=\overline{f_1(\overline{x})}.
\label{eq:conjuguaisons}\end{equation}
\label{proposition:poco}\end{Proposition}
\noindent {\bf DŽmonstration.} 
Notons que pour une fonction $X$ suffisamment dŽrivable~: \[\displaystyle{\frac{\partial}{\partial z_1}(X\circ{\cal
C})=\left(\frac{\partial X}{\partial z_2}\right)\circ{\cal C}},\quad
\displaystyle{\frac{\partial}{\partial z_2}(X\circ{\cal C})=\left(\frac{\partial X}{\partial z_1}\right)\circ{\cal C}}.\]
Soit $\rho\in{\cal T}$ une forme modulaire non nulle de poids $r$, symŽtrique ou
antisymŽtrique, posons
$\delta(\rho)=1$ si
$\rho$ est symŽtrique et $\delta(\rho)=-1$ si $\rho$ est antisymŽtrique. Alors~:
\[\displaystyle{\frac{\partial \rho}{\partial z_1}=\delta(r)\left(\frac{\partial \rho}{\partial z_2}\right)\circ{\cal C}}.\]
On en dŽduit que 
$[\varphi_2,\chi_5]_{1_1}\circ{\cal C}  = -[\varphi_2,\chi_5]_{1_2}$,
$[\varphi_2,\chi_6]_{1_1}\circ{\cal C}= [\varphi_2,\chi_6]_{1_2}$, et 
$[\chi_5,\chi_6]_{1_1}\circ{\cal
C} = -[\chi_5,\chi_6]_{1_2}$. 

Calculons maintenant $d_2(\rho)\circ{\cal C}$~:
\begin{eqnarray*}
\lefteqn{d_2(\rho)\circ{\cal C} =}\\
& = & [\rho,[\varphi_2,\chi_5]_{1_1}]_{1_2}\circ{\cal C}\\
& = & r(\rho\circ{\cal C})\left(\left(\frac{\partial }{\partial z_2}[\varphi_2,\chi_5]_{1_1}\right)\circ{\cal C}\right)-\\
& & 7 ([\varphi_2,\chi_5]_{1_1}\circ{\cal C})\left(\left(\frac{\partial \rho}{\partial z_2}\right)\circ{\cal C}\right)\\
& = & \delta(\rho)\left(r\rho\frac{\partial}{\partial z_1}([\varphi_2,\chi_5]_{1_1}\circ{\cal C})-
7(-[\varphi_2,\chi_5]_{1_2})\frac{\partial\rho}{\partial
z_1}\right)\\
& = & \delta(\rho)\left(-r\rho\frac{\partial}{\partial z_1}[\varphi_2,\chi_5]_{1_2}+7[\varphi_2,\chi_5]_{1_2}\frac{\partial\rho}{\partial
z_1}\right)\\
& = & -\delta(\rho)[\rho,[\varphi_2,\chi_5]_{1_2}]_{1_1}=-\delta(\rho)d_1(\rho).
\end{eqnarray*}
Une formule analogue est valide pour $d_2(\rho)\circ{\cal C}$. Donc si $s$ est une forme modulaire symŽtrique, alors $d^*(s)$ est
antisymŽtrique~: elle appartient donc ˆ $\chi_5{\cal T}_{\mbox{{\tiny sym}}}$. De mme, $d_*(s)$ est symŽtrique.
Si $a$ est une forme modulaire antisymŽtrique, alors $d^*(a)$ est symŽtrique. En raisonnant de la mme faon 
pour les autres dŽrivations $e^*,e_*,f^*,f_*$, on dŽmontre toutes les inclusions concernant ${\cal T}_{\mbox{{\tiny sym}}}$
et $\chi_5{\cal T}_{\mbox{{\tiny sym}}}$.

Les inclusions concernant ${\cal T}^*$ et $\chi_{15}{\cal T}^*$ en dŽcoulent aussi, car d'aprs 
la proposition \ref{proposition:structure2}, toute forme modulaire de Hilbert symŽtrique de poids parallles impairs
est produit d'une forme modulaire symŽtrique de poids pairs et de $\chi_{15}$. Comme $d^*,f^*$ envoient
des formes modulaires de poids pair sur des formes modulaires de poids impair, et $e^*$ envoie des formes 
modulaires de poids pair sur des formes modulaires de poids impair, on obtient 
toutes les images des sous-modules liŽs ˆ ${\cal T}_{\mbox{{\tiny sym}}}$ et ˆ ${\cal T}^*$. 

Ceci implique
que pour tout $r\in{\cal T}$ on a \[\overline{d^*(\overline{r})}=d^*(r),\quad\overline{d_*(\overline{r})}=-d_*(r).\]
On en dŽduit les ŽgalitŽs (\ref{eq:conjuguaisons}).
La dŽmonstration de la proposition \ref{proposition:poco} est terminŽe.

\begin{Corollaire} L'idŽal principal de ${\cal T}$ engendrŽ par $\chi_{15}$ est ${\mathfrak D}^*$-stable. L'idŽal engendrŽ par $\chi_{5}$ est $e^*$-stable.\label{corollaire:stable}\end{Corollaire}
\noindent {\bf DŽmonstration.} On utilise la collection d'inclusions (\ref{eq:ideal_principal_stable})~; on en dŽduit
les inclusions $d^*(\chi_{15}),e^*(\chi_{15}),f^*(\chi_{15})\in\chi_{15}{\cal T}^*$.

Si $e^*(\chi_5)$ est non nulle, alors elle a poids $15$, et est antisymetrique d'aprs
la proposition \ref{proposition:poco}. Comme toutes les formes modulaires de ${\cal T}$ de poids $\leq 14$
sont dans ${\cal T}^*$ (cf. proposition \ref{proposition:structure2}), on trouve
$e^*(\chi_5)\in\chi_5{\cal T}_0$. Le lemme qui suit implique que $e^*(\chi_5)\not=0$.

\begin{Lemme}
Si $t\in{\mathfrak D},y\in{\cal T}$ et $ty=0$, alors $y\in\CC$. Si $t^*\in{\mathfrak D},x\in{\cal T}$ et $t^*x=0$, alors $x\in\CC$.
\label{lemme:hortogonaux}\end{Lemme}
\noindent {\bf DŽmonstration.} 
Pour dŽmontrer ces propriŽtŽs, il suffit de le faire avec $x,y$ des formes modulaires, car toutes les 
dŽrivations sont isobares.

DŽmontrons la premire propriŽtŽ avec $t=d_1$~: la dŽmonstration donnŽe s'Žtend facilement 
aux autres dŽrivations de ${\mathfrak D}$.

Soit $y$ une forme modulaire telle que $d_1y=0$. On a $[y,[\varphi_2,\chi_5]_{1_2}]_{1_1}=0$. Ceci implique que $y$ et
$[\varphi_2,\chi_5]_{1_2}$ sont multiplicativement dŽpendants modulo $\CC^\times$ d'aprs le lemme 
11 de \cite{Pellarin:Hilbert}, 
mais ceci n'est possible que si $y\in\CC$, car
$y$ est  de poids parallles, alors que $[\varphi_2,\chi_5]_{1_2}$ ne l'est pas.

Soit maintenant $x\in{\cal T}$ une forme modulaire telle que $d^*x=0$.
On peut Žcrire $x=s+a$, o $s$ est une forme modulaire symŽtrique, et
$a$ est une forme modulaire antisymŽtrique. On a $d^*(s+a)=d^*s+d^*a=0$. D'aprs la proposition \ref{proposition:poco}, 
$d^*s$ est antisymŽtrique et
$d^*a$ est symŽtrique. Donc, on a $d^*s=0$ et $d^*a=0$. 

Il suffit de dŽmontrer que pour toute forme modulaire symŽtrique $s$, 
si $d^*s=0$ alors $s\in\CC$, et que pour toute
forme modulaire antisymŽtrique $a$, si $d^*a=0$ alors $a=0$.

Soit $s$ une forme modulaire symŽtrique non constante, telle que $d^*s=0$. Si $d_*s=0$, alors $d_1s=d_2s=0$, d'o $s\in\CC$,
d'aprs la premire partie du lemme. Donc $d_*s\not=0$.
On a $d_*s=d_1s=-d_2s$, car $d_1=d_*+d^*$ et $d_2=d_*-d^*$.
Donc, d'aprs (\ref{eq:conjuguaisons})~:
\begin{eqnarray*}
d_*s & = & d_1s=-d_2s\\
& = & -\iota d_1\iota s\\
& = & -\iota d_1s\\
& = & -\iota d_*s\\
& = & -d_*s,
\end{eqnarray*}
car $d_*$ preserve la symŽtrie (proposition \ref{proposition:poco}), d'o une contradiction. Donc $d^*s\not=0$. 

Soit $a$ une forme modulaire antisymŽtrique non constante, telle que $d^*a=0$. Si $d_*a=0$, alors $d_1a=d_2a=0$, d'o $a\in\CC$
et donc $a=0$ car elle est antisymŽtrique. Donc $d_*a\not=0$.
On a $d_*a=d_1a=-d_2a$ et donc, d'aprs (\ref{eq:conjuguaisons})~:
\begin{eqnarray*}
d_*a & = & d_1a=-d_2a\\
& = & -\iota d_1\iota a\\
& = & \iota d_1a\\
& = & \iota d_*a\\
& = & -d_*a,
\end{eqnarray*}
d'o une contradiction. Donc $d^*a\not=0$. La mme dŽmonstration fonctionne avec la dŽrivation $f^*$, et 
ˆ des petites modifications prs, avec $e^*$.

\section{IdŽaux diffŽrentiellement stables.}
Dans ce paragraphe nous dŽmontrons les thŽormes \ref{theoreme:ideaux_stables12} et \ref{theoreme:ideaux_stables}~: commenons par
le thŽorme \ref{theoreme:ideaux_stables}.

La dŽmonstration du thŽorme \ref{theoreme:ideaux_stables} procde en plusieurs Žtapes. Nous dŽ\-terminons (1) les idŽaux principaux stables pour ${\mathfrak D}^*$,
(2) les idŽaux premiers de ${\cal T}^*$ qui sont ${\mathfrak D}^*$-stables, n'ayant pas d'ŽlŽments
dans $\CC[\varphi_2]\cup\CC[\chi_5]\cup\CC[\chi_6]$,
(3) les idŽaux premiers de ${\cal T}^*$ qui sont ${\mathfrak D}^*$-stables, et qui contiennent au moins un ŽlŽment non nul de 
$\CC[\varphi_2]\cup\CC[\chi_5]\cup\CC[\chi_6]$.

\subsection{IdŽaux principaux ${\mathfrak D}^*$-stables.\label{section:principaux}}

\begin{Lemme}
L'idŽal $(\chi )$ est l'unique idŽal premier principal non nul de ${\cal T}^*$, qui est ${\mathfrak D}^*$-stable.
\label{lemme:ideaux_principaux}\end{Lemme}
\noindent {\bf DŽmonstration.} 
Le corollaire \ref{corollaire:stable} implique que $(\chi)$ est ${\mathfrak D}^*$-stable~: 
en effet, pour tout $t^*\in{\mathfrak D}^*$, $t^*(\chi)=t^*(\chi_{15}^2)=2\chi_{15}t^*(\chi_{15})$.
Or, $t^*(\chi_{15})\in{\chi_{15}}{\cal T}^*$ (corollaire \ref{corollaire:stable}), donc $t^*(\chi)\in\chi{\cal T}^*$.
C'est facile de dŽmontrer que $(\chi)$ est un idŽal premier de ${\cal T}^*$.

Nous dŽmontrons que $(\chi )$ est le seul idŽal premier principal $d^*$-stable de ${\cal T}^*$. 
Soit $p\in{\cal T}^*$ non constant tel que
$(d^*p)\subset(p)$. En particulier, il existe $a\in{\cal T}^*$ tel que $d^*p=ap$. 
On peut alors Žcrire~: $p=p_m+\cdots+p_l$ et $a=a_h+\cdots+a_k$ avec $p_i,a_i$
des formes modulaires de poids parallles $i$, 
$m\leq l$, $h\leq k$, $p_m,p_l\not=0$. Nous pouvons aussi supposer que $a_h,a_k\not=0$ d'aprs le lemme 
\ref{lemme:hortogonaux}.
On a~:
\[d^*p_m+\cdots+d^*p_l=a_hp_m+\cdots+a_kp_l.\] Donc $d^*p_m=a_hp_m$ et $d^*p_l=a_kp_l$. Comme $d^*$ est de poids $9$, 
on trouve que $h=k=9$, et $a=a_k=a_h$. Donc
$d^*p=ap$ et $a$ est une forme modulaire de poids parallles $9$. 
Mais $M_9(\Gamma)\cong\CC \varphi_2^2\chi_5$ (proposition \ref{proposition:structure2}), donc $a$ est un multiple non nul de $\varphi_2^2\chi_5$.

Ainsi, les deux formes modulaires $d^*\chi /\chi $ et $d^*p/p$ sont $\CC$-linŽairement dŽpendantes.
On a donc~:
\begin{equation}
\frac{d^*\chi }{\chi }=\lambda\frac{d^*p}{p},\label{eq:pour_instant_complexe}
\end{equation} avec $\lambda\in\CC^\times$. Nous montrons maintenant que $\lambda\in\QQ^\times$.

Jusqu'ˆ prŽsent, nous n'avons utilisŽ que des arguments algŽbriques~:
nous avons essentiellement travaillŽ avec des anneaux de polyn™mes. Il s'agit maintenant
d'utiliser la structure analytique des formes modulaires.

Pour achever notre dŽmonstration, 
il suffit de dŽmontrer que si $F,G,H$ sont des formes modulaires de Hilbert
algŽbriquement indŽpendantes, de
poids parallles respectivement $f,2,5$, avec $H$ parabolique et $G(\infty)=1$, alors
\[[F,G,H]^*=\lambda FG^2H,\]
implique $\lambda\in\QQ$, o nous avons posŽ~:
\begin{eqnarray*}
{[}F,G,H{]}^* & = & \frac{1}{2}({[}F,{[}G,H{]}_{1_2}{]}_{1_1}+{[}F,{[}G,H{]}_{1_1}{]}_{1_2}),
\end{eqnarray*}
car $d^*F=[F,\varphi_2,\chi_5]^*$ par dŽfinition.
Nous utilisons les sŽries de Fourier du lemme 15 de \cite{Pellarin:Hilbert}~: ici 
${\cal O}_{K,+}^*$ dŽsigne l'ensemble des ŽlŽments totalement positifs
du dual ${\cal O}_K^*$ de ${\cal O}_K$ pour la trace $\tr$ de $K$ sur $\QQ$. En particulier, en Žcrivant~:
\begin{eqnarray*}
F(z,z') & = & \sum_{\nu\in{\cal O}_{K,+}^*\cup\{0\}}f_\nu e(\nu z+\nu' z'),\\
G(z,z') & = & \sum_{\nu\in{\cal O}_{K,+}^*\cup\{0\}}g_\nu e(\nu z+\nu' z'),\\
H(z,z') & = & \sum_{\nu\in{\cal O}_{K,+}^*}h_\nu e(\nu z+\nu' z'),
\end{eqnarray*}
o $e(t)=e^{2\pi{\rm i}t}$, $(z,z')\in{\cal H}^2$, et avec $f_\nu,g_\nu,h_\nu\in\CC,g_0=1$, on trouve que
\begin{eqnarray}
[F,G,H]^* & = & \frac{1}{2}({[}F,{[}G,H{]}_{1_2}{]}_{1_1}+{[}F,{[}G,H{]}_{1_1}{]}_{1_2})\nonumber\\
& = & \sum_{\tau\in{\cal O}_{K,+}^*}e(\tau z+\tau' z')\sum_{\alpha+\nu+\mu=\tau}^\sharp
c_{\alpha,\nu,\mu}f_\alpha g_{\nu}h_{\mu},\label{eq:crochet_fourier}
\end{eqnarray}
o la somme $\sum^\sharp$ est indexŽe par des ŽlŽments $(\alpha,\nu,\mu)$ de $({\cal O}_{K,+}^*\cup\{0\})^3$, et o
\[c_{\alpha,\nu,\mu}:=\tr((f\alpha'-7(\nu'+\mu'))(2\nu-5\mu)).\]
De la mme faon, on calcule~:
\begin{equation}
FG^2H=\sum_{\tau\in{\cal O}_{K,+}^*}e(\tau z+\tau' z')\sum_{x+y+z+t=\tau}^\sharp
f_xg_yg_zh_t.\label{eq:produit_fourier}\end{equation}
Supposons pour commencer que $F$ soit parabolique.
Comme $F,H$ ne sont pas nulles, il existe deux entiers rationnels $k_0,h_0>0$ minimaux avec la propriŽtŽ
que pour quelques $\alpha\in{\cal O}_{K,+}^*$ avec $\tr(\alpha)=k_0$, on ait $f_\alpha\not=0$,
et pour quelques $\mu\in{\cal O}_{K,+}^*$ avec $\tr(\mu)=h_0$, on ait $h_\mu\not=0$.

Soient maintenant $\alpha_0,\mu_0\in{\cal O}_{K,+}^*$ tels que $\tr(\alpha_0)=k_0,\tr(\mu_0)=h_0$, avec la propriŽtŽ 
que $\alpha_0,\mu_0$ soient les plus grands possibles avec $f_{\alpha_0},h_{\mu_0}\not=0$ (On considre 
$K$ comme Žtant plongŽ dans $\RR$), soit $\tau_0=\alpha_0+\mu_0$.
L'expression (\ref{eq:produit_fourier}) et le choix de $\tau_0$ impliquent que le coefficient de Fourier de $FG^2H$ associŽ
ˆ $\tau_0$ est Žgal ˆ $g_0^2f_{\alpha_0}h_{\mu_0}=f_{\alpha_0}h_{\mu_0}\not=0$.

De mme, l'expression (\ref{eq:crochet_fourier}) et le choix de $\tau_0$ impliquent que le coefficient de Fourier de 
$[F,G,H]^*$ associŽ
ˆ $\tau_0$ est Žgal ˆ \begin{eqnarray*}c_{\alpha_0,0,\mu_0}g_0f_{\alpha_0}h_{\mu_0}& = & c_{\alpha_0,0,\mu_0}f_{\alpha_0}h_{\mu_0}\\
&=&\tr(5\mu_0(f\alpha_0'-7\mu_0))f_{\alpha_0}h_{\mu_0}.\end{eqnarray*}
Mais $c_{\alpha,\nu,\mu}$ Žtant la trace d'un ŽlŽment de $K$, c'est un nombre rationnel. Si $[F,G,H]^*=\lambda FG^2H$, alors
$c_{\alpha_0,0,\mu_0}f_{\alpha_0}h_{\alpha_0}=\lambda f_{\alpha_0}h_{\alpha_0}$, d'o $c_{\alpha_0,0,\mu_0}=\lambda\in\QQ$.

Supposons ensuite que $F$ ne soit pas parabolique~: on a alors $f_0\not=0$. Soit $\mu_0$ comme ci-dessus.
Le coefficient de Fourier associŽ ˆ $\mu_0$ dans (\ref{eq:produit_fourier}) est $f_0h_{\mu_0}$, tandis que
le coefficient de Fourier associŽ ˆ $\mu_0$ dans (\ref{eq:crochet_fourier}) est $f_0h_{\mu_0}c_{0,0,\mu_0}=
-70f_0h_{\mu_0}\no(\mu_0)$~: on en dŽduit que $\lambda\in\QQ^\times$.

Nous avons donc dŽmontrŽ l'existence de deux entiers rationnels $u,v$ tels que $d^*(p^u 
\chi ^v)/(p^u \chi ^v)=0$ (o $d^*$
designe  l'extension de la dŽrivation au corps des fractions de ${\cal T}^*$), 
ou de manire Žquivalente, $d^*(p^u \chi ^v)=0$. Le lemme \ref{lemme:hortogonaux} implique $p^u \chi ^v\in\CC$. On trouve
donc
$(p)=( \chi )$, et le lemme
\ref{lemme:ideaux_principaux} est dŽmontrŽ.

\medskip

\noindent {\bf Remarques}. Dans la dŽmonstration du lemme \ref{lemme:ideaux_principaux} nous
avons utilisŽ le fait remarquable
que
$\dim_\CC M_9(\Gamma)=1$. L'unicitŽ d'un idŽal principal $t^*$-stable 
pour une dŽriva\-tion $t^*$ est fausse en gŽnŽral. Par
exemple, l'opŽrateur $e^*$ possde l'idŽal principal premier stable $(\chi_5)$.

Observons aussi que $\chi_5,\chi_{15}$ sont des gŽnŽrateurs du ${\cal T}_0$-module libre ${\cal T}$. 
Dans \cite{Eichler:Projective} et \cite{Eichler:Acta},
on fait des hypothses trs gŽnŽrales quant'ˆ la structure des anneaux des formes modulaires. L'hypothse
p. 88 de \cite{Eichler:Acta} affirme que, pour tout anneau ${\mathfrak I}$ de formes modulaires
(dans une classe considŽrŽe par Baily et Borel dans \cite{Baily:Compactification}), il existe un certain sous-module
${\mathfrak H}$ tel que ${\mathfrak I}$ soit un ${\mathfrak H}$-module libre.
Ce serait trs intŽressant de comparer cette structure avec la structure diffŽrentielle
de ${\mathfrak I}$.

\subsection{IdŽaux premiers ${\mathfrak D}^*$-stables de hauteur $2$.\label{section:hauteur_2}}
Soient $F,G,H$ trois formes modulaires de Hilbert algŽbriquement indŽ\-pendantes, de poids parallles $f,g,h$.
D'aprs la proposition 3 de \cite{Pellarin:Hilbert}, la forme modulaire $[F,G,H] $, de poids parallles 
$f+g+h+2$ est non nulle. 

La formule suivante est ŽlŽmentaire et peut se dŽmontrer avec un calcul direct
(nous l'avons trouvŽe en suivant les techniques introduites dans \cite{Zagier:Indian}~; un logiciel de
calcul symbolique quelconque permet de la vŽrifier)~:
\begin{eqnarray}
\lefteqn{(f+g)(f+h){[}F,G,H{]} ^2 =}\label{eq:formulone}\\
& = &(f+g)(g+h)({[}F,F,G{]}^*{[}H,G,H{]}^*-{[}F,G,H{]}^*{[}H,F,G{]}^*)-\nonumber\\ & &
(g+h)(f+h)({[}H,G,H{]}^*{[}G,F,G{]}^*-{[}H,F,G{]}^*{[}G,F,H{]}^*)+\nonumber\\ & &
(g+h)^2({[}F,F,H{]}^*{[}G,G,H{]}^*-{[}F,G,H{]}^*{[}G,F,H{]}^*)\nonumber.
\end{eqnarray}
Soit $M$ la matrice~:
\[M=\left(\begin{array}{ccc} {[}G,G,H {]}^* & {[}H,G,H {]}^* & {[}F,G,H {]}^*\\ 
{[}G,G,F {]}^* & {[}H,G,F {]}^* & {[}F,G,F {]}^*\\ 
{[}G,H,F {]}^* & {[}H,H,F {]}^* & {[}F,H,F {]}^*\end{array}\right).\]
La formule (\ref{eq:formulone}) implique qu'elle a rang $\geq 2$. En effet, l'indŽpendance
algŽbrique de $F,G,H$ implique, en utilisant la proposition 3 de \cite{Pellarin:Hilbert}, que $[F,G,H] \not=0$,
donc $[F,G,H]^2\not=0$. Or, l'expression ˆ droite de (\ref{eq:formulone}) est une combinaison linŽaire de mineurs
de taille $2$ de la matrice $M$. Ainsi, un mineur de $M$ de taille $2$ est non nul 
(\footnote{On peut dŽmontrer que si $F,G,H$ sont des formes 
modulaires de poids parallle algŽbriquement indŽpendantes, alors $M$ est de rang $3$, mais nous ne le ferons pas 
ici.}).

Supposons que $F=\chi_6,G=\varphi_2$ et $H=\chi_5$. Alors
\[M=\left(\begin{array}{ccc} d^*\varphi_2 & d^*\chi_5 & d^*\chi_6\\ e^*\varphi_2 &
e^*\chi_5 & e^*\chi_6\\  f^*\varphi_2 & f^*\chi_5 & f^*\chi_6\end{array}\right).\] Soit $\tilde{M}$ la matrice 
des mineurs de $M$ de taille $2$.
D'aprs le raisonnement prŽcedant, cette matrice possde des coefficients non nuls.
La matrice suivante dŽcrit les poids des coefficients de $M$~:
\[\left(\begin{array}{ccc} 11 & 14 & 15\\ 12 & 15 & 16\\ 15 & 18 & 19\end{array}\right).\] Donc les coefficients de $\tilde{M}$ sont
tous des formes modulaires de Hilbert, et les poids sont repartis de la manire suivante~:
\[\left(\begin{array}{ccc} 34 & 31 & 30\\ 33 & 30 & 29\\ 30 & 27 & 26\end{array}\right)\]
(on notera que tous les coefficients de l'adjointe $\hat{M}$ de $M$ ont mme poids $45$). Nous avons le lemme suivant.
\begin{Lemme}
Soit ${\cal P}$ un idŽal premier non principal de ${\cal T}^*$ tel que ${\cal
P}\cap(\CC[\varphi_2]\cup\CC[\chi_5]\cup\CC[\chi_6])=(0)$. 
Si ${\cal P}$ est ${\mathfrak D}^*$-stable, alors tous les coefficients de $\tilde{M}$ appartiennent ˆ ${\cal P}$. De plus, 
$ \chi \in{\cal P}$.
\label{lemme:non_principal}\end{Lemme}
\noindent {\bf DŽmonstration.} Les hypothses impliquent qu'il existe dans ${\cal P}$, trois poly\-n™mes irrŽductibles $A,B,C$ avec les 
propriŽtŽs suivantes~:
\begin{eqnarray*}
\frac{\partial A}{\partial \chi_6} = 0, & & \frac{\partial A}{\partial \varphi_2}\frac{\partial A}{\partial \chi_5} \not= 0 \\
\frac{\partial B}{\partial \chi_5} = 0, & & \frac{\partial B}{\partial \varphi_2}\frac{\partial B}{\partial \chi_6} \not= 0 \\
\frac{\partial C}{\partial \varphi_2} = 0, & & \frac{\partial C}{\partial \chi_5}\frac{\partial C}{\partial \chi_6} \not= 0.
\end{eqnarray*}
On peut mme supposer que les degrŽs partiels de ces polyn™mes soient les plus petits possibles,
car il n'y a pas d'idŽal non trivial de ${\cal T}^*$ qui soit stable pour les trois dŽrivŽes partielles $\partial/\partial \varphi_2,
\partial/\partial \chi_5,\partial/\partial \chi_6$.

On a $d^*A,d^*B,\ldots,f^*C\in{\cal P}$. En particulier, $d^*A,e^*A\in{\cal P}$, et \[(d^*A)(e^*\varphi_2)-(e^*A)(d^*\varphi_2)\in{\cal
P}.\] On vŽrifie l'ŽgalitŽ~:
\[(d^*A)(e^*\varphi_2)-(e^*A)(d^*\varphi_2)=\frac{\partial A}{\partial \chi_5}(d^*\chi_5e^*\varphi_2-e^*\chi_5d^*\varphi_2).\] 
Comme $\partial A/\partial \chi_5\not\in{\cal P}$,
on a $d^*\chi_5e^*\varphi_2-e^*\chi_5d^*\varphi_2\in{\cal P}$. 
De la mme faon que ci-dessus, on vŽrifie alors que tous les coefficients de la matrice $\tilde{M}=(\tilde{m}_{i,j})_{i,j}$
appartiennent ˆ ${\cal P}$. 

La dernire propriŽtŽ du lemme se dŽmontre en utilisant la formule (\ref{eq:formulone}), car
$[\chi_6,\varphi_2,\chi_5] ^2$ est proportionnel ˆ $\chi $ (avec constante de proportionnalitŽ non nulle), et Žgal ˆ une combinaison 
linŽaire non triviale des cofficients qui se trouvent sur l'anti-diagonale de $\tilde{M}$.

\medskip

\noindent {\bf Remarque.}
Les arguments de ce paragraphe pourraient se gŽnŽraliser
au cas o $K$ est un autre corps quadratique rŽel diffŽrent de $\QQ(\sqrt{5})$.

\subsection{DŽtermination explicite de rŽlations diffŽrentielles.}

Pour continuer, nous devons expliciter les rŽlations
engendrŽes par les dŽrivations de ${\mathfrak D}^*$. Nous utilisons la proposition suivante.
\begin{Proposition}
On a les rŽlations suivantes.
\begin{eqnarray*}
d^*(\varphi_{2}) & = & \frac{4}{5}\chi_5(\varphi_2^3-1050\chi_6)\\
d^*(\chi_6) & = & 
-\frac{2}{5}\chi_5(\varphi_2^2\chi_6+875\chi_5^2)\\ 
d^*(\chi_5) & = & 
\frac{1}{10}\varphi_2(7\chi_6^2-15\varphi_2 \chi_5^2)
\end{eqnarray*}
\begin{eqnarray*} 
e^*(\varphi_2) & = & -1152\chi_6^2+240\varphi_2\chi_5^2+\frac{4}{5}\varphi_2^3\chi_6\\
e^*(\chi_6) & = & -240\chi_5^2\chi_6-\frac{8}{5}\varphi_2^2\chi_6^2+\frac{4}{5}\varphi_2^3\chi_5^2\\
e^*(\chi_5) & = & \chi_5\left(-\frac{6}{5}\varphi_2^2\chi_6+200\chi_5^2\right)
\end{eqnarray*}
\begin{eqnarray*}
f^*(\varphi_2) & = & \chi_5 \left(550 \chi_5^2 - \frac{4}{5}\varphi_2^2 \chi_6\right)\\
f^*(\chi_6) & = & \chi_5 \left(\frac{11}{4}\varphi_2^2 \chi_5^2 - \frac{59}{20}\varphi_2 \chi_6^2\right)\\
f^*(\chi_5) & = & \chi_6 \left(\frac{7}{2}\varphi_2 \chi_5^2 - \frac{33}{10} \chi_6^2\right).
\end{eqnarray*}
\label{proposition:systeme_differentiel}\end{Proposition}
\noindent {\bf Esquisse de dŽmonstration}. 
Montrons comment obtenir la premire des rŽlations de la proposition \ref{proposition:systeme_differentiel}.

Nous savons, d'aprs le lemme \ref{lemme:hortogonaux}, que $d^*(\varphi_2)$ est non nulle. Il s'agit d'une forme 
modulaire de poids parallle $11$. La proposition \ref{proposition:structure2} implique que~:
\[d^*(\varphi_2)=\alpha\chi_5\varphi_2^3+\beta\chi_5\chi_6,\quad\alpha,\beta\in\CC.\]
Cette relation linŽaire implique l'existence de relations linŽaires entre coefficients de Fourier
de $d^*(\varphi_2),\chi_5\varphi_2^3,\chi_5\chi_6$, et le calcul explicite de quelques coefficient
de Fourier est suffisant pour calculer $\alpha,\beta$.
Pour les autres relations aussi, on effectue une rŽsolution numŽrique de systmes linŽaires, en utilisant la proposition
\ref{proposition:structure2}. 

Plue en gŽnŽral, il faut calculer un certain nombre de coefficients de Fourier
de $\varphi_2,\chi_5,\chi_6,\chi_{15}$.
Pour ce faire il suffit, d'aprs les relations diffŽren\-tielles de la 
proposition \ref{proposition:structure2}, de conna"tre tous les coefficients de Fourier ˆ l'infini de $\varphi_2,\chi_5$, correspondants ˆ des ŽlŽments
totalement positifs du dual pour la trace de l'anneau des entiers de $K$, de trace majorŽe par un entier assez
grand (nous les appellons coefficients de \og petite trace \fg). 

La connaissance de ces coefficients de Fourier de $\varphi_2$ nous permet de dŽterminer une partie 
de la sŽrie de Fourier de $\chi_6$ et de $\chi_5^2$. 
Pour calculer les coefficients de Fourier de petite trace de $\varphi_2$, on peut appliquer les formules de
\cite{Pellarin:Hilbert}.

Pour calculer les coefficients de Fourier de petite trace
de $\chi_5$, on peut appliquer la formule explicite de \cite{Bruinier:Borcherds}. Tous ces coefficients
permettent ˆ leur tour de calculer, via le crochet $[\cdot,\cdot,\cdot] $, les coefficients de 
Fourier de petite trace pour $\chi_{15}$.

La mŽthode est la mme que dans \cite{Pellarin:Hilbert} ou \cite{Resnikoff:Automorphic}, nous ne donnons 
pas les dŽtails de calcul ici~: le logiciel \og Pari \fg  permet 
de vŽrifier facilement toutes les rŽlations
de la proposition, et nous rŽmercions D. Simon pour avoir vŽrifiŽ toutes ces formules. 

\medskip

\noindent {\bf Remarque.} Nous pouvons expliciter les rŽlations du corollaire \ref{corollaire:stable}
gr‰ce ˆ la proposition \ref{proposition:systeme_differentiel}~:
\begin{eqnarray*}
d^* \chi  & = & -2\varphi_2^2\chi_5 \chi \\
e^* \chi  & = & -4(\varphi_2^2\chi_6-300\chi_5^2) \chi \\
f^* \chi  & = & -\varphi_2\chi_5\chi_6 \chi .
\end{eqnarray*}

Dans toute la suite, nous allons utiliser les rŽlations 
dŽcrites explicitement par la proposition \ref{proposition:systeme_differentiel}.

\subsection{IdŽaux premiers ${\mathfrak D}^*$-stables de hauteur $3$.}
Dans ce paragraphe, nous Žtudions les idŽaux de ${\cal T}^*$ qui ont intersection non nulle avec
l'ensemble $\CC[\varphi_2]\cup\CC[\chi_5]\cup\CC[\chi_6]$~; nous verrons qu'ils ont tous hauteur $\geq 3$. 

\begin{Lemme}
Soit ${\cal P}$ un idŽal premier de ${\cal T}^*$,
supposons qu'au moins une des conditions suivantes soit vŽrifiŽe~:
\begin{equation}
{\cal P}\cap\CC[\varphi_2]\not=(0),\quad{\cal P}\cap\CC[\chi_5]\not=(0),\quad{\cal
P}\cap\CC[\chi_6]\not=(0).\label{eq:condition_1}\end{equation}
Si ${\cal P}$ est ${\mathfrak D}^*$-stable, alors ${\cal P}$ contient l'idŽal premier $(\chi_5,\chi_6)$.
En particulier, $\chi\in{\cal P}$.
\label{lemme:intersections}\end{Lemme}
\noindent {\bf DŽmonstration.} 
Comme ${\cal P}$ est premier, on a ou bien $\varphi_2-c_2\in{\cal P}$, ou bien $\chi_5-c_5\in{\cal P}$, ou bien $\chi_6-c_6\in{\cal P}$,
pour $c_2,c_5,c_6\in\CC$. 

\medskip

\noindent {\bf (1)}.
Supposons d'abord que $\chi_5-c_5\in{\cal P}$~: nous montrons que $\chi_6\in{\cal P}$. On a~:
\begin{equation}
d^*\chi_5,f^*\chi_5,d^*{}^2\chi_5,f^*{}^2\chi_5\in{\cal P},\label{eq:d5}\end{equation}
car pour tout $t^*\in{\mathfrak D}^*$, $t^*(\CC)=0$.
En appliquant la proposition \ref{proposition:systeme_differentiel}, on peut calculer explicitement les polyn™mes
(\ref{eq:d5}). 

Pour vŽrifier les calculs suivants, il est conseillŽ de se servir d'un logiciel de calcul symbolique.
On calcule les rŽsultants~:
\begin{eqnarray}
\mbox{RŽs}{}_{\chi_5}(d^*\chi_5,d^*{}^2\chi_5) & = & \xi_1\varphi_2^5\chi_6^6A_1\in{\cal P}\label{eq:resu51}\\
\mbox{RŽs}{}_{\chi_5}(f^*\chi_5,f^*{}^2\chi_5) & = & \xi_2\chi_6^{15}A_2\in{\cal P},\label{eq:resu52}
\end{eqnarray}
o $\xi_1,\xi_2\in\QQ^\times$, et o
$A_1,A_2\in\CC[\varphi_2,\chi_6]$ dŽsignent des polyn™mes isobares non nuls de poids $12$. On vŽrifie que $A_1,A_2$ sont premiers entre eux,
et pour des 
raisons de poids, ils sont Žgaux ˆ des produits de la forme
\begin{equation}\prod_s(a_s\varphi_2^3-\chi_6),\label{eq:prod6}\end{equation}
o l'on vŽrifie que $a_s\in\bar{\QQ}^\times$ pour tout $s$.

Comme ${\cal P}$ est premier, dans chacun des deux polyn™mes 
de (\ref{eq:resu51}) et (\ref{eq:resu52}) il existe au moins un facteur irrŽductible
appartenant ˆ ${\cal P}$~; si $\chi_6\in{\cal P}$ nous avons 
terminŽ, supposons donc par l'absurde que $\chi_6\not\in{\cal P}$.

Ainsi, ${\cal P}$ contient au moins un facteur de $A_2$, d'aprs (\ref{eq:resu52}), et ne peut pas contenir $\varphi_2$
d'aprs (\ref{eq:resu51}), car sinon, en Žliminant $\varphi_2$ d'un facteur
de la forme (\ref{eq:prod6}), nous obtenons $\chi_6\in{\cal P}$. Donc ${\cal P}$ contient $A_1,A_2$, et comme il est premier, il contient aussi 
deux facteurs irrŽductibles premiers entre eux
de la forme (\ref{eq:prod6}). En Žliminant $\varphi_2$, nous parvenons une fois de plus ˆ une contradiction, car nous obtenons $\chi_6\in{\cal P}$.

\medskip

\noindent {\bf (2)}. Supposons que $\chi_6-c_6\in{\cal P}$~: nous montrons que $\chi_5\in{\cal P}$. On a
que~:
\[d^*\chi_6,d^*{}^2\chi_6,e^*f^*\chi_6\in{\cal P}.\] On calcule cette fois les rŽsultants suivants, qui appartiennent tous ˆ ${\cal P}$~:
\begin{eqnarray*}
\mbox{RŽs} {}_{\chi_6}(d^*\chi_6,d^*{}^2\chi_6) & = & \xi_3\varphi_2^3\chi_5^7B_1\in{\cal P}\label{eq:resu61}\\
\mbox{RŽs} {}_{\chi_6}(d^*\chi_6,e^*f^*\chi_6) & = & \xi_4\chi_5^9B_2\in{\cal P},\label{eq:resu62}
\end{eqnarray*}
o $\xi_3,\xi_4\in\QQ^\times$ et $B_1,B_2$ sont deux polyn™mes isobares non nuls de $\CC[\varphi_2,\chi_5]$ de poids $10,20$.

On vŽrifie que $B_1,B_2$ sont premiers entre eux~; pour de raisons de poids, ce sont des produits de la forme~:
\begin{equation}\prod_s(b_s\varphi_2^5-\chi_5^2),\label{eq:prod5}\end{equation}
avec $b_s\in\bar{\QQ}^\times$ pour tout $s$.
En suivant la technique appliquŽe au point {\bf (1)}, on obtient $\chi_5\in{\cal P}$.

\medskip

Dans les deux cas {\bf (1)} ou {\bf (2)} on voit que l'idŽal $(\chi_5,\chi_6)$ est contenu dans ${\cal P}$. En effet, d'aprs {\bf (1)}, 
si $\chi_5\in{\cal P}$, alors 
$\chi_6\in{\cal P}$ et $(\chi_5,\chi_6)\subset{\cal P}$. Si parcontre $\chi_5-c_5\in{\cal P}$ et $c_5\not=0$, alors on trouve $\chi_6\in{\cal P}$ d'aprs {\bf (1)},
puis $\chi_5\in{\cal P}$ d'aprs {\bf (2)}, et finalement, ${\cal P}={\cal T}^*\supset(\chi_5,\chi_6)$.

\medskip

\noindent {\bf (3)}. Supposons pour terminer que $\varphi_2-c_2\in{\cal P}$. On calcule les rŽsultants~:
\begin{eqnarray*}
\mbox{RŽs} {}_{\varphi_2}(d^*\varphi_2,d^*{}^2\varphi_2) & = & \xi_5\chi_5^{11}C_1\in{\cal P}\label{eq:resu21}\\
\mbox{RŽs} {}_{\varphi_2}(e^*\varphi_2,e^*{}^2\varphi_2) & = & \xi_6\chi_5^8C_2\in{\cal P},\label{eq:resu22}
\end{eqnarray*}
avec $\xi_5,\xi_6\in\QQ^\times$, et $C_1,C_2$ deux polyn™mes isobares non nuls de poids $30,60$, premiers entre eux de la forme
\[\prod_s(d_s\chi_5^6-\chi_6^5),\]
avec $d_s\in\bar{\QQ}^\times$ pour tout $s$.
En utilisant la mme technique des deux points prŽcedants, on voit que $\chi_5\in{\cal P}$, d'o $(\chi_5,\chi_6)\subset{\cal P}$
gr‰ce au point {\bf (1)}.

\medskip

Dans tous les cas, $(\chi_5,\chi_6)\subset{\cal P}$. Comme $\chi\in(\chi_5,\chi_6)$ (formule (\ref{eq:klein})), 
le lemme \ref{lemme:intersections} est entirement
dŽmontrŽ.

\subsection{Preuve des thŽormes \ref{theoreme:modulaires}, \ref{theoreme:ideaux_stables12} et \ref{theoreme:ideaux_stables}.}
\noindent {\bf DŽmonstration du thŽorme \ref{theoreme:ideaux_stables}}. Soit ${\cal P}$ un idŽal premier non trivial de ${\cal T}^*$, ${\mathfrak D}^*$-stable.
Sa hauteur gŽomŽtrique $h$ varie dans l'ensemble $\{1,2,3\}$. Si $h=1$ alors ${\cal P}$ est principal, donc Žgal ˆ $(\chi )$
d'aprs le lemme \ref{lemme:ideaux_principaux}. Supposons maintenant que $h\geq 2$. Si ${\cal P}
\cap(\CC[\varphi_2]\cup\CC[\chi_5]\cup\CC[\chi_6])\not=(0)$,
d'aprs le lemme \ref{lemme:intersections} on a $ \chi \in{\cal P}$, et la proposition est dŽmontrŽe dans ce cas.

Nous pouvons donc supposer que ${\cal P}\cap(\CC[\varphi_2]\cup\CC[\chi_5]\cup\CC[\chi_6])=(0)$ et que ${\cal P}$ est non principal, et
d'aprs le lemme \ref{lemme:non_principal}, $ \chi \in{\cal P}$~: la preuve du thŽorme \ref{theoreme:ideaux_stables}
est terminŽe.

\medskip

\noindent {\bf DŽmonstration du thŽorme \ref{theoreme:ideaux_stables12}.} Nous devons d'abord Žnoncer trois lemmes techniques.
\begin{Lemme}
Soit ${\cal Q}$ un idŽal non nul de ${\cal T}$. Si ${\cal Q}\cap\overline{{\cal Q}}=(0)$, alors ${\cal Q}$ est principal.
\label{lemme:int_nulle_id_principal}\end{Lemme}
\noindent {\bf DŽmonstration.} Soit ${\cal Q}$ un idŽal non principal de ${\cal T}$. Il existe deux ŽlŽments non nuls $p_1,p_2$ de 
${\cal Q}$ tels que pour tout choix d'ŽlŽments $q_1,q_2\in{\cal T}$ non tous nuls, alors $q_1p_1+q_2p_2\not=0$. Supposons par l'absurde que 
${\cal Q}\cap\overline{{\cal Q}}=(0)$~: on a en particulier ${\cal Q}\cap{\cal T}^*=(0)$, car sinon, si $s\in{\cal Q}\cap{\cal T}$
est un ŽlŽment non nul, alors $\overline{s}=s$ et $s\in\overline{{\cal Q}}$.

Donc $p_1,p_2\not\in{\cal T}^*$~: nous pouvons Žcrire alors
$p_1=x_1+\chi_{15}y_1,p_2=x_2+\chi_{15}y_2$, avec $x_1,x_2,y_1,y_2\in{\cal T}^*$, et $y_1,y_2$ non nuls. 
Mais $y_2p_1-y_1p_2=y_2x_1-y_1x_2\not=0$ est un ŽlŽment de
${\cal Q}$ qui est dans ${\cal T}^*$, d'o une contradiction.
\begin{Lemme}
Soit ${\cal Q}$ un idŽal premier non nul de ${\cal T}$, ne contenant pas $\chi_{15}$, soit ${\cal P}={\cal Q}\cap{\cal T}^*$. On a~:
\begin{equation}
{\cal Q}\cap\overline{{\cal Q}}={\cal P}{\cal T}.\label{eq:intersection_ideal_premier}\end{equation}
\label{lemme:intersection_ideal_premier}\end{Lemme}
\noindent {\bf DŽmonstration.} 
C'est clair que ${\cal Q}\cap\overline{{\cal Q}}\supset{\cal P}{\cal T}$~; montrons que
${\cal Q}\cap\overline{{\cal Q}}\subset{\cal P}{\cal T}$. Clairement, ${\cal Q}\cap{\cal T}^*=\overline{{\cal Q}}\cap{\cal T}^*$,
car $\overline{{\cal T}^*}={\cal T}^*$. D'autre part, ${\cal P}$ est un idŽal premier de ${\cal T}^*$.

Soit $x\in{\cal Q}\cap\overline{{\cal Q}}$. On a $\overline{x}\in{\cal Q}$ et donc \[\beta:=x+\overline{ x} \in{\cal P}.\]

D'autre part, aussi $x-\overline{x}\in{\cal Q}$~: or $x-\overline{x}\in \chi_{15}{\cal T}^*$ et il existe $\alpha\in{\cal T}^*$ tel que
$x-\overline{x}=\chi_{15}\alpha$~; cet ŽlŽment appartient ˆ ${\cal Q}$. 

Nous avons $\chi_{15}^2,\alpha\in{\cal T}^*$, donc $\chi_{15}^2\alpha\in{\cal P}$ Or, par hypothse, $\chi_{15}^2\not\in{\cal Q}$,
donc $\chi_{15}^2\not\in{\cal P}$ car ${\cal Q}$ est premier. Finalement, $\alpha\in{\cal P}$, d'o~:
\[\chi_{15}\alpha=x-\overline{x}\in{\cal P}{\cal T}.\] 

On a donc que $\beta$ et $\chi_{15}\alpha$ appartiennent ˆ ${\cal P}{\cal T}$. 
Ainsi, \[x=(\beta+\chi_{15}\alpha)/2\in{\cal P}{\cal T},\] d'o ${\cal P}{\cal T}\supset{\cal
Q}\cap\overline{{\cal Q}}$, d'o l'ŽgalitŽ
(\ref{eq:intersection_ideal_premier})~: le lemme \ref{lemme:intersection_ideal_premier} est dŽmontrŽ.

\begin{Lemme} Il existe des nombres $l_1,l_2,l_3\in K^\times$ tels que 
si $r\in{\cal T}^*$, alors~:
\begin{equation}d_*r=l_1\chi_{15}\frac{\partial r}{\partial \chi_6},\quad e_*r =
l_2\chi_{15}\frac{\partial r}{\partial \chi_5},\quad f_*r = l_3\chi_{15}\frac{\partial
r}{\partial
\chi_2}.\label{eq:derivees_partielles}\end{equation} En particulier on a~:
\begin{equation}
d_*\chi_{15}=\frac{l_1}{2}\frac{\partial \chi }{\partial \chi_6},\quad
e_*\chi_{15}= \frac{l_2}{2}\frac{\partial \chi }{\partial \chi_5},\quad
f_*\chi_{15}= \frac{l_3}{2}\frac{\partial \chi }{\partial \chi_2}.\label{eq:derivees_partielles_chi15}
\end{equation}
\label{lemme:chi15_deri}\end{Lemme}
\noindent {\bf DŽmonstration.} Utilisons 
la proposition 3 de \cite{Pellarin:Hilbert}. On a $d_*(\varphi_2)=d_*(\chi_5)=0$, et
$d_*(\chi_6)=(7/2)[\chi_6,\varphi_2,\chi_5] \not=0$, car $\chi_6,\varphi_2,\chi_5$ sont algŽbriquement in\-dŽpendantes. 
Donc $d_*(\chi_6)$ est proportionnel ˆ
$\chi_{15}$ pour une constante de proportionnalitŽ non nulle~: plus prŽ\-cisement, on a $d_*(\chi_6)=(11/\sqrt{5})\chi_{15}$.
Les autres formules de (\ref{eq:derivees_partielles}) s'obtiennent de la mme faon.
Les constantes $l_1,l_2,l_3$ peuvent se calculer explicitement.

Les formules explicites de (\ref{eq:derivees_partielles_chi15}), s'obtiennent 
tout simplement en calculant $d_*(\chi ),e_*(\chi ),f_*(\chi )$, car on a, par exemple pour la dŽrivation $d_*$~:
$d_*(\chi_{15}^2)=2\chi_{15}d_*(\chi_{15})$. 

\medskip

\noindent {\bf Preuve du thŽorme \ref{theoreme:ideaux_stables12}.} Soit ${\cal Q}$ un idŽal premier de ${\cal T}$, 
stable pour tous les opŽrateurs de ${\mathfrak D}$. On a   que $\overline{{\cal Q}}$ est ${\mathfrak D}$-stable, d'aprs (\ref{eq:conjuguaisons}). 
En particulier, les idŽaux ${\cal Q}$ et $\overline{{\cal Q}}$ sont ${\mathfrak D}^*$-stables.
Comme ${\mathfrak D}^*{\cal T}^*\subset{\cal T}^*$, l'idŽal 
${\cal P}={\cal Q}\cap{\cal T}^*$ est ${\mathfrak D}^*$-stable. Nous distinguons maintenant deux cas,
suivant que ${\cal Q}\cap\overline{{\cal Q}}=(0)$ ou non.

\medskip

\noindent {\bf (1)}. Supposons que ${\cal Q}\cap\overline{{\cal Q}}=(0)$. Alors, d'aprs le lemme \ref{lemme:int_nulle_id_principal}, ${\cal Q}$ est
un idŽal principal~: soit $p$ un gŽnŽrateur~: alors $p\not\in{\cal T}^*$. 
Il existe six formes modulaires $a_1,a_2,b_1,b_2,c_1,c_2$, de poids parallles respectivement $9,9,10,10,13,13$, telles
que~:
\[d_ip=a_ip,\quad e_ip=b_ip,\quad f_ip=c_ip\mbox{ pour }i=1,2.\]
En particulier, $a_i,b_i,c_i\in{\cal T}^*$ ($i=1,2$) et satisfont $\overline{a}_i=a_i,\ldots$,
car toutes les formes modulaires de ${\cal T}$ de poids $\leq 14$ se trouvent dans ${\cal T}^*$, d'aprs la proposition
\ref{proposition:structure2}. D'aprs les ŽgalitŽs (\ref{eq:conjuguaisons})~:
\begin{eqnarray*}
\overline{d_1p} & = & d_2(\overline{p})=a_1\overline{p},\\
\overline{d_2p} & = & d_1(\overline{p})=a_2\overline{p},
\end{eqnarray*}
Posons $\theta=p\cdot \overline{p}\in{\cal T}^*$. On a~:
\begin{eqnarray*}
d_1\theta & = & \overline{p}d_1(p)+pd_1(\overline{p})\\
& = & a_1\theta+a_2\theta\\
& = & (a_1+a_2)\theta\quad\mbox{ et de manire analogue~: }\\
d_2\theta & = & (a_2+a_1)\theta.
\end{eqnarray*}
Donc, $d^*\theta=a\theta$, avec $a=a_2+a_1$, et on obtient, de la mme faon~:
$e^*\theta=b\theta,f^*\theta=c\theta$, avec $b=b_1+b_2$ et $c=c_1+c_2$. Donc, l'idŽal $(\theta)$
de ${\cal T}^*$ est 
${\mathfrak D}^*$-stable, Žgal ˆ $(\chi ^l)$ pour un certain entier $l$, d'aprs le lemme \ref{lemme:ideaux_principaux}. 
Soit $p'$ un ŽlŽment de ${\cal T}$ 
tel que $p'\overline{p'}= \chi $. Alors $p'$ est une forme modulaire, et proportionnel ˆ $\chi_{15}$. Donc $l=1$ et $\theta\in\CC^\times\chi $.

Nous avons dŽmontrŽ que si ${\cal Q}$ est un idŽal principal ${\mathfrak D}$-stable, alors ${\cal Q}=(\chi_{15})$.
En utilisant l'ŽgalitŽ (\ref{eq:derivees_partielles_chi15}) du lemme \ref{lemme:chi15_deri}, on vŽrifie que
que $(\chi_{15})$ n'est pas ${\mathfrak D}_*$-stable~: comme il est quand-mme 
${\mathfrak D}^*$-stable, on en dŽduit que $(\chi_{15})$ n'est pas ${\mathfrak D}$-stable. Nous avons dŽmontrŽ en
fait~:
\begin{Lemme}
Il n'existe pas d'idŽal principal non nul ${\mathfrak D}$-stable dans ${\cal T}$.
\label{lemme:ideaux_princ_R}\end{Lemme}
\noindent {\bf (2).} Supposons que ${\cal Q}\cap\overline{{\cal Q}}\not=(0)$. 
D'aprs le lemme \ref{lemme:intersection_ideal_premier}, ${\cal P}$ est
non nul, et premier~; il est de plus ${\mathfrak D}^*$-stable. 
Supposons par l'absurde que $\chi_{15}\not\in{\cal Q}$. 

Comme ${\cal P}$ est ${\mathfrak D}^*$-stable, il contient $ \chi $ d'aprs le thŽorme 
\ref{theoreme:ideaux_stables}. Donc $ \chi \in{\cal Q}$ et $\chi_{15}\in{\cal Q}$~: une contradiction.
La dŽmonstration du thŽorme \ref{theoreme:ideaux_stables12} est terminŽe.

\medskip

\noindent {\bf Preuve du thŽorme \ref{theoreme:modulaires}}. Soit ${\mathfrak P}$ un idŽal premier stable
de ${\cal Y}(K)$, contenant une forme modulaire de Hilbert non nulle. Alors l'idŽal ${\cal Q}$
engendrŽ par toutes les formes modulaires de ${\mathfrak P}$ est 
un idŽal premier
non nul et ${\mathfrak D}$-stable de ${\cal T}$, qui contient $\chi_{15}$ d'aprs le
thŽorme \ref{theoreme:ideaux_stables12}.

\section{Un rŽsultat plus prŽcis.\label{section:paragrapheB}}
Soient $a,b$ deux nombres complexes,
considŽrons l'idŽal \[{\cal P}(a,b)=(a\varphi_2^5-\chi_5^2,b\varphi_2^3-\chi_6)\] de ${\cal T}^*$
et l'idŽal \[{\cal Q}(a,b)=({\cal P}(a,b),\chi_{15})\] de ${\cal T}$.
Ici nous donnons une esquisse de dŽmonstration du rŽsultat suivant.
\begin{Proposition} Si ${\cal P}$ est un idŽal premier non nul de ${\cal T}^*$ qui est ${\mathfrak D}^*$-stable
et de hauteur $\geq 2$, alors ${\cal P}$ contient ${\cal P}(a,b)$ avec~:
\[(a,b)\in{\cal E}:=\left\{\left(\frac{1}{800000},\frac{1}{800}\right),\left(\frac{1}{253125},\frac{1}{675}\right),(0,0)\right\}.\]
Si ${\cal Q}$ est un idŽal premier non nul de ${\cal T}$ qui est ${\mathfrak D}$-stable, alors il contient ${\cal Q}(a,b)$ avec
$(a,b)\in{\cal E}$.
\label{theoreme:ideaux_stables_ht_2}\end{Proposition}
La proposition \ref{theoreme:ideaux_stables_ht_2} est certes plus prŽcise des
thŽormes \ref{theoreme:ideaux_stables12} et \ref{theoreme:ideaux_stables}, mais la dŽmon\-stration
que nous en donnons utilise de manire essentielle la proposition \ref{proposition:systeme_differentiel}, ce qui 
ne peut que limiter des Žventuelles gŽnŽralisations.

Nous avons besoin du lemme qui suit, dont la dŽmonstration est seulement esquissŽe,
et dŽpend fortement des formules explicites de la proposition \ref{proposition:systeme_differentiel}.
\begin{Lemme} Nous avons les propriŽtŽs suivantes.
\begin{enumerate}
\item Si $a,b\in\CC^\times$, ou si $a=b=0$, l'idŽal ${\cal P}(a,b)$ est primaire de hauteur $2$. 
\item L'idŽal ${\cal P}(a,b)$ est ${\mathfrak D}^*$-stable si et seulement si $(a,b)\in{\cal E}$.
\item Si $(a,b)\not\in{\cal E}$, alors l'idŽal engendrŽ par ${\cal P}(a,b)$ et ${\mathfrak D}^*{\cal P}(a,b)$ a hauteur $\geq 3$.
\end{enumerate}
\label{lemme:idealPab}\end{Lemme}
\noindent {\bf Esquisse de dŽmonstration.} {\bf (1)}. On vŽrifie que le radical de ${\cal P}(a,b)$ dans ${\cal T}^*$ pour 
$a,b\in\CC^\times$ est l'idŽal premier de hauteur $2$~:
\[(b\varphi_2^3-\chi_6,a\varphi_2^2\chi_6-b\chi_5^2,b^2\varphi_2\chi_5^2-a\chi_6^2).\]
Si $a=b=0$, le radical de ${\cal P}(0,0)$ est l'idŽal $(\chi_5,\chi_6)$.

\medskip 

\noindent {\bf (2)}. Nous voulons trouver tous les couples $(a,b)\in\CC^2$ tels que 
\[d^*(a\varphi_2^5-\chi_5^2),d^*(b\varphi_2^3-\chi_6)\in{\cal P}(a,b).\]
On obtient, ˆ l'aide de la proposition \ref{proposition:systeme_differentiel}~:
\begin{eqnarray*}
d^*(a\varphi_2^5-\chi_5^2) & = & -3\varphi_2^2\chi_6(a\varphi_2^5-\chi_5^2)+\\
& & \frac{7}{5}\varphi_2\chi_6(\chi_6 + \varphi_2^3(3000  a + b))(b\varphi_2^3-\chi_6)-\\
& & \frac{7}{5}\varphi_2^7\chi_6(b^2 + 3000  ab-5a),\\
d^*(b\varphi_2^3-\chi_6) & = & -350\chi_5(a\varphi_2^5-\chi_5^2)+\\
& & \frac{2}{5}\varphi_2^2\chi_5(6300 b-1)(b\varphi_2^3-\chi_6)+\\
& & \frac{14}{5}\varphi_2^5\chi_5(125 a + b - 900 b^2).
\end{eqnarray*}
Pour que $d^*(a\varphi_2^5-\chi_5^2),d^*(b\varphi_2^3-\chi_6)\in{\cal P}(a,b)$, il faut
et il suffit que $125 a + b - 900 b^2=b^2 + 3000  ab-5a=0$, c'est-ˆ-dire, que $(a,b)$
appartienne ˆ ${\cal E}$. Sur la base de ces informations, on vŽrifie ensuite
que les idŽaux ${\cal P}(0,0),{\cal P}(1/800000,\\
1/800),{\cal P}(1/253125,1/675)$ sont
aussi ${\mathfrak D}^*$-stables, et que les idŽaux ${\cal Q}(0,0),\\
{\cal Q}(1/800000,1/800),{\cal Q}(1/253125,1/675)$ sont
${\mathfrak D}$-stables. 

\medskip

\noindent {\bf (3)}. On utilise, comme pour la propriŽtŽ {\bf (2)}, la proposition
\ref{proposition:systeme_differentiel}. On dŽmontre en calculant des rŽsultants, que si ${\cal P}(a,b)$
n'est pas ${\mathfrak D}^*$-stable, alors l'idŽal engendrŽ par ${\cal P}(a,b)$ et ${\mathfrak D}^*{\cal P}(a,b)$ 
contient les trois ŽlŽments $\varphi_2,\chi_5,\chi_6$.

\noindent {\bf DŽmonstration de la proposition \ref{proposition:systeme_differentiel}.} 
On vŽrifie que si $(a,b)\in{\cal E}$, alors ${\cal Q}(a,b)$ a hauteur $2$ car
$\chi_{15}^2\in{\cal P}(a,b)$.

Soit ${\cal P}$ un idŽal premier non nul de hauteur $\geq 2$, ${\mathfrak D}^*$-stable.
Suivant la dŽmonstration du lemme \ref{lemme:non_principal}, tous les coefficients $\tilde{m}_{i,j}$ de la matrice
$\tilde{M}$ appartiennent ˆ ${\cal P}$. 

La proposition \ref{proposition:systeme_differentiel}
permet de calculer $\tilde{m}_{i,j}$ explicitement. 
Pour $(i,j),(h,k)\in\{1,2,3\}$, $(i,j)\not=(h,k)$, on calcule ensuite le rŽsultant~:
\[R_{i,j,h,k,\varrho}:=\mbox{RŽs}_{\varrho}(\tilde{m}_{i,j},\tilde{m}_{h,k}),\quad\varrho=\chi_5,\chi_6,\]
qui sont des ŽlŽments isobares de ${\cal P}$~; on vŽrifie, en utilisant de prŽference un logiciel de calcul symbolique, que
ces ŽlŽments sont tous non nuls.

On peut factoriser explicitement ces rŽsultants dans $\bar{\QQ}$. Pour des raisons de poids, on trouve que, pour tout
$i,j,h,k$ comme ci-dessus~:
\begin{eqnarray}
R_{i,j,h,k,\chi_6} & = & \lambda\varphi_2^\alpha\chi_5^\beta\prod_{s=1}^t(a_s\varphi_2^5-\chi_5^2)\label{eq:factorisation}\\
R_{i,j,h,k,\chi_5} & = & \mu\varphi_2^{\alpha'}\chi_6^{\gamma'}\prod_{s=1}^{t'}(b_s\varphi_2^3-\chi_6),\label{eq:factorisation2}
\end{eqnarray}
avec $\alpha,\beta,\ldots,t'$ entiers positifs, $\lambda,\mu\in\CC^\times$, 
$a_s,b_s$ des nombres complexes non nuls,
le tout dŽpendant de $i,j,k,h$. 

Les rŽsultants $R_{i,j,h,k,\varrho}$
appartiennent tous ˆ ${\cal P}$, et comme ${\cal P}$ est premier, au moins un facteur de $R_{i,j,h,k,\chi_6}$ ˆ droite de (\ref{eq:factorisation})
et au moins un facteur de $R_{i,j,h,k,\chi_5}$ ˆ droite de (\ref{eq:factorisation2}), se trouvent dans ${\cal P}$.

D'aprs le lemme \ref{lemme:intersections}, nous pouvons supposer que $\varphi_2,\chi_5,\chi_6\not\in{\cal P}$, car dans ce cas,
nous savons dŽjˆ que
${\cal P}(0,0)\subset{\cal P}$.

Donc ${\cal P}$ contient un idŽal de la forme ${\cal P}(a,b)$ avec $a,b$ complexes non nuls. Si la hauteur de ${\cal
P}$ est $3$, alors ${\cal P}\cap(\CC[\varphi_2]\cup\CC[\chi_5]\cup\CC[\chi_6])\not=\{0\}$, et
${\cal P}(0,0)=(\chi_5^2,\chi_6)\subset{\cal P}$ (lemme \ref{lemme:intersections}).
Supposons que
${\cal P}$ ait hauteur
$2$. L'idŽal ${\cal P}(a,b)$ ne pouvant pas tre principal, il doit avoir hauteur $2$ et doit tre ${\mathfrak D}^*$-stable,
d'aprs le lemme \ref{lemme:idealPab}.

Le lemme \ref{lemme:idealPab} implique $(a,b)\in{\cal E}$, ce qui prouve la partie de l'enoncŽ de la proposition 
\ref{theoreme:ideaux_stables_ht_2} concernant l'idŽal ${\cal P}$.

La propriŽtŽ concernant les idŽaux ${\mathfrak D}$-stables de ${\cal T}$, se dŽduit de ce que nous venons de dŽmontrer,
car d'aprs le lemme \ref{lemme:int_nulle_id_principal}, ${\cal T}$ n'a pas d'idŽal non nul principal et
${\mathfrak D}$-stable.

\section{Appendice.}
\begin{Lemme} Soit $F$ une forme quasi-modulaire de Hilbert
de poids $(k_1,\ldots,k_n)$ et profondeur $s$. Soit 
\begin{equation}P=
\sum_{s_1+\cdots+s_n\leq s} f_{s_1,\ldots,s_n}(z) X_1^{s_1}\cdots X_n^{s_n}\label{eq:associated}\end{equation}
le polyn™me $P$ associŽ ˆ $F$ dans (\ref{eq:definition_hilbert}).
Alors $k\geq 2s$, et pour tout $s_1,\ldots,s_n$ tel que $s_1+\cdots+s_n=s$,
la fonction $f_{s_1,\ldots,s_n}(z)$ est une forme modulaire 
de Hilbert de poids 
$(k_1-2s_1,\ldots,k_n-2s_n)$.
\label{lemma:structure}\end{Lemme}
\noindent {\bf DŽmonstration.} Soient $A,B$ deux ŽlŽments de $\Gamma_K$ et Žcrivons~:
\[A=\sqm{a}{b}{c}{d},\quad B=\sqm{\alpha}{\beta}{\gamma}{\delta},\quad AB=\sqm{u}{v}{x}{y},\]
soit $F$ comme dans les hypothses du lemme.
Nous avons~:
\begin{eqnarray*}
\lefteqn{F(AB(z))=}\\
& = & \prod_{i=1}^n(c_iB_i(z_i)+d_i)^{k_i}
\sum_{s_1+\cdots+s_n\leq s}c_{s_1,\ldots,s_n}(B(z))\times\\ & &\prod_{i=1}^n\left((\gamma_i
z_i+\delta_i)^2\left(\frac{x_i}{x_iz_i+y_i}-\frac{\gamma_i}{\gamma_iz_i+\delta_i}\right)\right)^{s_i},
\end{eqnarray*}
donc:
\begin{eqnarray*}
\lefteqn{F(AB(z))=}\\
& = & \prod_{i=1}^n(c_iB_i(z_i)+d_i)^{k_i}\sum_{s_1+\cdots+s_n\leq s}c_{s_1,\ldots,s_n}(B(z))\times\\
& & \prod_{i=1}^n(\gamma_i
z_i+\delta_i)^{2s_i}\times\\
& &
\sum_{t_i=0}^{s_i}\binomial{s_i}{t_i}(-1)^{s_i-t_i}\left(\frac{x_i}{x_iz_i+y_i}\right)^{t_i}\left(
\frac{\gamma_i}{\gamma_iz_i+\delta_i}\right)^{s_i-t_i}
\end{eqnarray*}
\begin{eqnarray*}
& = &
\prod_{i=1}^n(c_iB_i(z_i)+d_i)^{k_i}\times\\
& &\sum_{s_1=0}^s\sum_{s_2=0}^{s-s_1}\cdots\sum_{s_n=0}^{s-s_1-\cdots-s_{n-1}}
c_{s_1,\ldots,s_n}(B(z))\prod_{i=1}^n(\gamma_i z_i+\delta_i)^{2s_i}\times\\
& &
\sum_{t_i=0}^{s_i}\binomial{s_i}{t_i}(-1)^{s_i-t_i}\left(\frac{x_i}{x_iz_i+y_i}\right)^{t_i}\left(
\frac{\gamma_i}{\gamma_iz_i+\delta_i}\right)^{s_i-t_i}
\end{eqnarray*}
\begin{eqnarray*}
& = & \prod_{i=1}^n(c_iB_i(z_i)+d_i)^{k_i}\sum_{t_1+\cdots+t_n\leq s}\prod_{i=1}^n
\left(\frac{x_i}{x_iz_i+y_i}\right)^{t_i}\times\\
& & \sum_{s_1=t_1}^s\sum_{s_2=t_2}^{s-s_1}\cdots\sum_{s_n=t_n}^{s-s_1-\cdots-s_{n-1}}
f_{s_1,\ldots,s_n}(B(z))\prod_{i=1}^n(\gamma_i z_i+\delta_i)^{2s_i}\times\\
& &\binomial{s_i}{t_i}(-1)^{s_i-t_i}\left(\frac{\gamma_i}{\gamma_iz_i+\delta_i}\right)^{s_i-t_i}.
\end{eqnarray*}
D'autre part:
\begin{eqnarray*}
\lefteqn{F(AB(z))=}\\
 &= & \prod_{i=1}^n(x_iz_i+y_i)^{k_i}
\sum_{s_1+\cdots+s_n\leq s}f_{s_1,\ldots,s_n}(z)\prod_{i=1}^n\left(\frac{x_i}{x_iz_i+y_i}\right)^{s_i}.\end{eqnarray*}
Ces deux expressions diffŽrentes de $F(AB(z))$ doivent co•ncider pour tout $A,B$ et $z$.
Nous pouvons fixer pour un instant $B$, $z$ en laissant varier $A$ de telle sorte que l'ensemble~:
\[\left\{\left(\frac{q_1}{q_1z_1+r_1},\ldots,\frac{q_n}{q_nz_n+r_n}\right)\right\}\] soit Zariski dense
dans ${\cal H}^n$. Comme
\[\prod_{i=1}^n(x_iz_i+y_i)^{k_i}=\prod_{i=1}^n((c_iB_i(z_i)+d_i)(\gamma_iz_i+\delta_i))^{k_i},\]
nous obtenons l'idŽntitŽ formelle:
\begin{eqnarray*}
\lefteqn{\prod_{i=1}^n(\gamma_iz_i+\delta_i)^{k_i}
\sum_{t_1+\cdots+t_n\leq s}f_{t_1,\ldots,t_n}(z)X_1^{t_1}\cdots X_n^{t_n}=}\\
& = & \sum_{t_1+\cdots+t_n\leq s}X_1^{t_1}\cdots X_n^{t_n}\times\\
& & \sum_{s_1=t_1}^s\sum_{s_2=t_2}^{s-s_1}\cdots\sum_{s_n=t_n}^{s-s_1-\cdots-s_{n-1}}
f_{s_1,\ldots,s_n}(B(z))\prod_{i=1}^n(\gamma_i z_i+\delta_i)^{2s_i}\times\\
& &\binomial{s_i}{t_i}(-1)^{s_i-t_i}\left(\frac{\gamma_i}{\gamma_iz_i+\delta_i}\right)^{s_i-t_i}\\
& = & \sum_{t_1+\cdots+t_n\leq s}X_1^{t_1}\cdots X_n^{t_n}\times\\
& & \sum_{s_1=t_1}^s\sum_{s_2=t_2}^{s-s_1}\cdots\sum_{s_n=t_n}^{s-s_1-\cdots-s_{n-1}}
f_{s_1,\ldots,s_n}(B(z))\times\\
& &\prod_{i=1}^n\binomial{s_i}{t_i}(-\gamma_i)^{s_i-t_i}\left(\gamma_iz_i+\delta_i\right)^{s_i+t_i}.
\end{eqnarray*}
En comparant les coefficients des mon™mes $X_1,\ldots,X_n$ nous obtenons, pour tout
$(t_1,\ldots,t_n)$ tel que $t_1+\cdots+t_n\leq s$:
\begin{eqnarray}
f_{t_1,\ldots,t_n}(z)& = &\sum_{s_1=t_1}^s\sum_{s_2=t_2}^{s-s_1}\cdots\sum_{s_n=t_n}^{s-s_1-\cdots-s_{n-1}}
f_{s_1,\ldots,s_n}(B(z))\times\label{eq:formuletta}\\
& &\prod_{i=1}^n\binomial{s_i}{t_i}(-\gamma_i)^{s_i-t_i}(\gamma_iz_i+\delta_i)^{s_i+t_i-k_i}
\nonumber\end{eqnarray}
Si $(t_1,\ldots,t_n)$ est tel que $t_1+\cdots+t_n=s$ et $f_{t_1,\ldots,t_n}(z)\not=0$, alors
la somme indŽxŽe par les $s_1,\ldots,s_n$ ˆ droite de l'ŽgalitŽ (\ref{eq:formuletta}) n'a qu'un seul terme,
avec $s_i=t_i$, $i=1,\ldots,n$. Nous obtenons, pour ce terme~:
\[f_{s_1,\ldots,s_n}(B(z))=\prod_{i=1}^n(\gamma_iz_i+\delta_i)^{k_i-2s_i}f_{s_1,\ldots,s_n}(z),\]
ŽgalitŽ valide pour tout $B,z$, ce qui implique que
$f_{s_1,\ldots,s_n}(z)$ est une forme modulaire de Hilbert de poids
$(k_1-2s_1,\ldots,k_n-2s_n)$~; la dŽmonstration du lemme est complte.

\medskip

Comme consŽquence de ce qui precde, on voit que si $F$ est une forme quasi-modulaire de Hilbert non constante, de
poids $(k_1,\ldots,k_n)$, alors $k_i\geq 0$ pour tout $i=1,\ldots,n$. De plus, une forme quasi-modulaire de Hilbert
de poids $(0,\ldots,0)$ est une fonction constante.
\begin{Lemme}
Supposons que $[K:\QQ]>1$ et soit $F$ une forme quasi-modulaire de Hilbert dont le poids $(k_1,\ldots,k_n)$ est non nul
et satisfait
$\prod_{i=1}^nk_i=0$~: alors $F=0$.
\label{lemma:non_existence}\end{Lemme}
\noindent {\bf DŽmonstration.} Supposons par l'absurde
qu'il existe une forme quasi-mo\-dulaire de Hilbert $F$ non nulle, de poids
$(k_1,\ldots,k_n)$ et de profondeur $s$, avec la propriŽtŽ que $k_j>0$ et $k_h=0$~; 
soit $P$ le polyn™me associŽ (\ref{eq:associated}). Nous allons procŽder en deux Žtapes.

\medskip

\noindent {\bf (1)}.
Nous dŽmontrons que $k_i$ est pair pour tout $i=1,\ldots,n$, que
$\sum_ik_i=2s$, et qu'il n'existe qu'un seul $n$-uple d'entiers positifs $(s_1,\ldots,s_n)$ avec
$\sum_is_i=s$ et $f_{\und{s}}\not=0$ dans (\ref{eq:associated}). Nous prouvons de plus que $f_{\und{s}}$ est constante.

Clairement, nous pouvons aussi supposer que $s>0$, car sinon $F$ est une forme modulaire de Hilbert non nulle,
en contradiction avec le lemme 6.3 p. 18 de \cite{Geer:Hilbert}.

Le lemme \ref{lemma:structure} nous dit que $\sum_ik_i\geq 2s$ et que
pour tout $(s_1,\ldots,s_n)$ tel que $\sum_is_i=s$ et $f_{\und{s}}\not=0$, alors $k_i\geq 2s_i$. 

Si $s$ satisfait
$\sum_ik_i> 2s$, alors pour tout $(s_1,\ldots,s_n)$ tel que $\sum_is_i=s$ et $f_{\und{s}}\not=0$, 
il existe $j$ tel que $k_j>2s_j$. De plus, par hypothse, il existe $h$
tel que $k_h=0$, ce qui implique
$s_h=0$. Donc $k_j-2s_j>0$ et $k_h-2s_h=0$.

Le lemme \ref{lemma:structure} implique que $f_{s_1,\ldots,s_n}$ est une forme 
modulaire de Hilbert de poids $(k_1-2s_1,\ldots,k_n-2s_n)$ et le lemme 6.3 p. 18 de \cite{Geer:Hilbert} implique 
$f_{s_1,\ldots,s_n}=0$. Mais alors, dans ce cas, le degrŽ de $P$ est $<s$ d'o une contradiction.

La seule possibilitŽ est donc que $\sum_ik_i=2s$, $k_i\geq 2s_i$ (ce qui implique $k_i=2s_i$
pour $i=1,\ldots,n$ et $k_i$ pair). Il n'existe qu'un $n$-uple d'entiers positifs $(s_1,\ldots,s_n)$ tel que
$\sum_is_i=s$, et $f_{\und{s}}\not=0$. D'aprs le lemme (\ref{lemma:structure}), la forme modulaire de Hilbert $f_{s_1,\ldots,s_n}$ est 
de poids $\und{0}$ et c'est donc
une constante non nulle $\lambda\in\CC$.

\medskip

\noindent {\bf (2)}. Etant donnŽ un poids $(k_1,\ldots,k_n)$ tel que
$k_j>0$ et $k_h=0$, nous dŽmontrons que l'espace vectoriel ${\cal Y}(K)_{(k_1,\ldots,k_n)}$ a dimension au plus $1$ sur $\CC$.

\medskip

En effet, si $F,G\in{\cal Y}(K)_{(k_1,\ldots,k_n)}$ sont non nulles, alors elles ont profondeur $(1/2)\sum_ik_i$ 
et \[f_{k_1/2,\ldots,k_n/2}(z)=\lambda,\quad g_{k_1/2,\ldots,k_n/2}(z)=\mu,\] avec $\lambda,\mu\in\CC^\times$.
Donc~:
\[\mu F-\lambda G\in{\cal Y}(K)_{(k_1,\ldots,k_n)}\]
est une forme quasi-modulaire de Hilbert de profondeur $<\sum_{i}k_i$, qui est nulle d'aprs
l'Žtape (1).

\medskip

Nous terminons la dŽmonstration du lemme~: soit $F$ non nul comme dans les hypothses. 
Le poids de $F^3$ et le poids de $D_1^{k_1}\cdots D_n^{k_n}F$ sont Žgaux ˆ $(3k_1,\ldots,3k_n)$~:
ces fonctions sont linŽairement indŽpendantes, et il existe un nombre complexe $\tau$ tel que l'Žquation
aux dŽrivŽes partielles suivante soit satisfaite~:
\[D_1^{k_1}\cdots D_n^{k_n}F=\tau F^3.\] Mais $F$ admet une expansion en sŽrie de Fourier ˆ l'infini~;  en comparant
les coefficients de cette expansion on trouve $F=0$, d'o une contradiction qui termine la preuve du lemme.

\medskip

DŽmontrons le thŽorme \ref{theorem:theorem1}. Soit ${\mathfrak E}$ le
sous-ensemble de
$\NN^n$ dont les points $(k_1,\ldots,k_n)$ sont tels que ${\cal Y}(K)_{(k_1,\ldots,k_n)}\not=(0)$.

Supposons par l'absurde que ${\cal Y}(K)$ soit de type fini, engendrŽ par des ŽlŽment
$F_1,\ldots,F_m$. On peut supposer que $F_1,\ldots,F_m$ soient des formes quasi-modulaires de Hilbert
de poids $\und{k}_1,\ldots,\und{k}_m$. Alors~:
\[{\mathfrak  E}=\NN\und{k}_1+\cdots+\NN\und{k}_m.\]
Soit ${\mathfrak  C}$ l'enveloppe convexe de ${\mathfrak  E}$ de $\RR^n$~: alors ${\mathfrak  C}=\RR_{\geq 0}\und{k}_1+\cdots+\RR_{\geq
0}\und{k}_m$ est un c™ne polyhedral qui ne contient aucun angle plan 
d'amplitude $\geq
\pi/2$, car d'aprs le lemme \ref{lemma:non_existence}, si $(k_1,\ldots,k_n)\in{\mathfrak  E}$, alors $\prod_ik_i\not=0$.

Donc ${\mathfrak  C}$ n'est pas stable par les translations par les vecteurs de la base canonique de
$\RR^n$. 
Mais ${\cal Y}(K)$ est un anneau diffŽrentiel~; donc ${\mathfrak  E}$ doit tre stable par tous les ŽlŽmŽnts de $\NN^n$,
et ${\mathfrak  C}$ Žgalement. Ceci mne ˆ une contradiction, et la preuve du thŽorme 
\ref{theorem:theorem1} est terminŽe.

\end{document}